\documentclass[12pt,thmsa,emstex]{article}
\usepackage{amsfonts}
\usepackage[utf8]{inputenc}
\usepackage[T1]{fontenc}
\usepackage{t1enc}
\usepackage{amssymb,bbm}
\def\1{\mathbbm{1}}
\usepackage{epsfig}
\usepackage{amsmath}
\usepackage[french,english]{babel}
\usepackage{amsmath}
\usepackage{graphics}
\usepackage{layout}
\usepackage{latexsym}
\usepackage{color}
\usepackage{verbatim}
\usepackage{soul}

\def\Z{\mathbb{Z}}

\def\mylabel#1{\label{#1}}

\newtheorem{theorem}{Theorem}[section]
\newtheorem{lemma}[theorem]{Lemma}
\newtheorem{corollary}[theorem]{Corollary}
\newtheorem{proposition}[theorem]{Proposition}

\newtheorem{definition}[theorem]{Definition}
\newtheorem{exercise}[theorem]{Exercise}
\newtheorem{hypothesis}{Hypothesis}[section]

\newtheorem{remark}[theorem]{Remark}
\newtheorem{remarks}{Remarks}[section]
\newtheorem{example}[theorem]{Example}

\def\bit{\begin{itemize}}
\def\eit{\end{itemize}}

\reversemarginpar   

\def\bc{\begin{center}}
\def\ec{\end{center}}

\def\bthm{\begin{theorem}}
\def\ethm{\end{theorem}}

\def\bcor{\begin{corollary}}
\def\ecor{\end{corollary}}

\def\bprop{\begin{proposition}}
\def\eprop{\end{proposition}}

\def\blem{\begin{lemma}}
\def\elem{\end{lemma}}

\def\bex{\begin{example}}
\def\eex{\end{example}}

\def\bexo{\begin{exercise}}
\def\eexo{\end{exercise} }

\def\brem{\begin{remark}}
\def\erem{\end{remark}}
\def\brems{\begin{remarks}}
\def\erems{\end{remarks}}
\def\prf{{\bf Proof: }}

\def\bdes{\begin{description}}
\def\edes{\end{description}}

\def\ita{\item[(a)]}
\def\itb{\item[(b)]}

\def\itc{\item[(c)]}

\def\iti{\item[(i)]}
\def\itii{\item[(ii)]}

\def\itiii{\item[(iii)]}

\def\beq{\begin{equation}}
\def\eeq{\end{equation}}

\def\ben{\begin{enumerate}}
\def\een{\end{enumerate}}

\def\beqar{\begin{eqnarray}}
\def\eeqar{\end{eqnarray}}

\def\beqarr{\begin{eqnarray*}}
\def\eeqarr{\end{eqnarray*}}

\def\qed{\hfill $\Box$ \\[2ex]}
\def\prf{{\bf Proof: }\hspace{.1in}}

\def\cM{\mathcal{M}}

\catcode`\@=11

\newcommand{\B}{\mathcal{B}}
\newcommand{\C}{\mathcal{C}}

\newcommand{\cD}{{\mathcal D}}

\newcommand{\bD}{{\overline{\cD}}}

\def\Ind{{\mathbbm{1}}}
\def\ZZ{{\mathbb Z}}       
\def\RR{{\mathbb R}}  
\def\Rp{{\mathbb R}_+}   

\def\NN{{\mathbb N}}

\def\rar{\rightarrow}
\def\eps{\varepsilon}

\def\bla{}

\begin{document}

\title{Degenerate processes killed at the boundary of a domain}
\author{M Bena\"im$^1$, N Champagnat$^{2}$, W O\c cafrain$^{2}$,  D Villemonais$^{2,3}$}
\footnotetext[1]{Institut de Math\'{e}matiques, Universit\'{e} de Neuch\^{a}tel, Switzerland.}
\footnotetext[2]{Université de Lorraine, CNRS, Inria, IECL, F-54000 Nancy, France.}
\footnotetext[3]{Institut Universitaire de  France.}
\maketitle
\abstract{We investigate quasi-stationarity properties of  Feller processes that are killed
  when exiting a relatively compact set. Our main result provides general conditions ensuring that such a process possesses a
  (possibly non unique) quasi stationary distribution. Conditions ensuring uniqueness and exponential convergence are discussed.
Our conditions are well-suited to the study of degenerate processes, such as nonelliptic diffusions or piecewise deterministic
    Markov processes (PDMP). The results are applied to  stochastic differential equations
and we illustrate the application to PDMPs with an example}.



\section{Introduction}
\label{sec:intro}

\subsection{Context and objectives}

This paper  investigates certain properties of a  Feller Markov process $(X_t)_{t \geq 0}$  living on some metric space $M,$ killed
when it exits an open, relatively compact set $D$, at time $\tau_D^{out}=\inf\{t\geq 0:X_t\not\in D\}$.

We discuss  conditions ensuring:
\bdes
\ita Existence of a {\em quasi-stationary distribution} (QSD). That is, a probability $\mu$ on $D$ such that
 $$\mathbb{P}_{\mu} (X_t \in \cdot \mid \tau_D^{out} > t) = \mu(\cdot);$$
  \itb Uniqueness of a QSD;
   \itc  Convergence  of the conditional laws  $\mathbb{P}_{\nu} (X_t \in \cdot \mid\tau_D^{out} > t)$for any
     probability measure $\nu$ on $D$ toward this QSD,  in total variation.
\edes
We refer the reader to the survey paper \cite{Meleard-Villemonais-2012} and the monograph \cite{collet2012quasi} for a comprehensive introduction to the  theory of quasi-stationarity.

Probabilistic criteria (based on Lyapunov functions and minorization conditions) implying $(c)$ are given in
\cite{champagnat2016exponential,champagnat2017general}. These apply to elliptic diffusions, as shown in \cite{champagnat2018criteria}
and in \cite[Section 4]{champagnat2017general}. Based on the well-known link between QSDs and eigenproblems for the adjoint
  Dirichlet operator (see e.g.~\cite{Meleard-Villemonais-2012}), quasi-stationarity has also been
widely studied through spectral-theoretic (in a wide sense) arguments (Sturm-Liouville theory~\cite{cattiaux2009quasi,littin2012uniqueness,kolb2012quasilimiting,HeningKolb2014}, compactness
  or quasi-compactness criteria~\cite{Gosselin2001,HinrichsKolbEtAl2020,GNW2020}, Hilbert-Schmidt
  theory~\cite{ChampagnatDiaconisEtAl2012}, Krein-Rutman theorem~\cite{FerreRoussetStoltz2021,collet2012quasi}, Tychonov's fixed
  point theorem~\cite{ColMartMel}, $R$-positive matrices~\cite{FerrariKestenEtAl1996}, orthogonal
  polynomials~\cite{KarlinMcGregor1957}, to cite few references among many). For diffusion processes (assuming that $D$
    is connected), a lot of works made use of
  minorization-Lyapunov
  methods~\cite{champagnat2017uniform,ChampagnatVillemonais2015a,DelMoralVillemonais2018,champagnat2018criteria,Velleret2018} or
  spectral methods
\cite{Pinsky1985,GongQianEtAl1988,cattiaux2009quasi,hening2021quasi,kolb2012quasilimiting,littin2012uniqueness,martinez1994quasi}.

However, the treatment of quasi-stationarity beyond elliptic diffusions i.e.\ non-elliptic diffusions is a more delicate
problem. This topic received only few contributions~\cite{GNW2020,LRR2021}, based on spectral
arguments or two-sided estimates. Concerning the conditions in
\cite{champagnat2016exponential,champagnat2017general}, they rely on certain global  parabolic Harnack type inequalities (condition (A2) in
\cite{champagnat2016exponential} or (E3) in \cite{champagnat2017general}) that are difficult to verify in general for degenerate
(non-elliptic) diffusions, although some results are known in particular cases~\cite{GolseImbertEtAl2019}.  The two
  references~\cite{GNW2020,LRR2021} study the question~$(c)$ for strong Feller diffusions.
As we will see, the strong Feller assumption, although natural in certain situations, is unnecessarily limiting. It is interesting to understand under what milder
conditions and for which processes more general than diffusions the existence property $(a)$ (respectively the uniqueness property $(b)$) holds
while $(b)$ (respectively $(c)$) does not, a situation which is typical of {\color{black} degenerate processes, i.e. processes which are Feller but not strong Feller.}

\bla{Other degenerate processes  in the previous sense include piecewise deterministic Markov processes (PDMP), which received contributions on
  quasi-stationary distributions only for
  specific examples of processes~\cite{champagnat2016exponential,BansayeCloezEtAl2019,VillemonaisWatson2022,CloezFritsch2022}.  The main
  objective of our paper is to propose a general framework and to show how it can be used
  to obtain new results on general degenerate diffusions, and to illustrate how in can be used for some PDMPs.}

The remainder of this introduction is devoted to the presentation of our main results in the specific case where $(X_t)$ is a
degenerate  diffusion process. Then, Section \ref{sec:main} sets up the general framework and addresses
points (a) (Theorem \ref{th:main}), (b) (Theorem \ref{th:uniqueness}) and (c) (Theorem \ref{th:convergence}). \bla{These results all
  assume that the Green kernel sends bounded continuous functions on continuous functions vanishing at the boundary of the domain.
  Existence (a) is proved under mild accessibility conditions; uniqueness (b) assuming irreducibility and the existence of a
  positive eigenfunction for the Green operator (a property that can be deduced for example from the compactness of the Green operator
  and Krein-Rutman's theorem); and convergence (c) assuming the existence of a positive eigenfunction for the Green operator and a
  small set condition. These results rely on Tychonov's fixed point theorem for existence and on a version of Harris theorem
  based on Lyapunov-minorization criteria applied to the $Q$-process (the process conditioned to never be absorbed, cf.\
  e.g.~\cite{Meleard-Villemonais-2012}).} These results are applied in Section \ref{sec:sde} to prove (among other things) the
results stated  in Section~\ref{sec:sde0} below for diffusions. The proofs rely on criteria for accessibility, and in particular on Stroock and
  Varadhan's support theorem, on Feller's property and the regularity of point of the boundary to prove that the Green kernel sends
  bounded continuous functions on continuous functions vanishing at the boundary, on Rotschild and Stein's H\"older
  estimates~\cite{rothschild1976} for solutions to hypoelliptic PDEs and on the existence of continuous solutions to
    the Dirichlet
  problem established by Bony~\cite{Bony} to prove that the Green kernel is compact, and on estimates on the Green kernel
  from~\cite{Bony} and on the transition densities by Ichihara and Kunita~\cite{Ichihara} to prove small-set-like properties. Section
\ref{sec:PDMP}  discusses a simple PDMP example which opens perspectives for other applications of our results.
\subsection{Description of the results for SDEs}
\label{sec:sde0}

{\color{black}Let $\cD \subset \mathbb R^n$ be an open connected set with compact closure $\bD$ and boundary $\partial \cD = \bD \setminus \cD.$}

Consider  the stochastic differential equation on $\mathbb R^n$
\beq
\label{eq:sde}
dX_t = S^0(X_t) dt + \sum_{j = 1}^m S^j(X_t) \circ dB_t^j,
 \eeq where $\circ$ refers to the Stratonovich stochastic integral, $S^0, S^j, j = 1, \ldots, m$ are smooth vector fields on $\mathbb
 R^n$ and $B^1, \ldots, B^m$ $m$ independent Brownian motions. 
 Equivalently, for the reader familiar with Ito's calculus,
 $$dX_t = \left[S^0(x) + \frac{1}{2} \sum_{j = 1}^m DS^j(x) S^j(x)\right] dt + \sum_{j = 1}^m S^j(X_t)  dB_t^j,$$ where
 $DS^j(x)$ stands for the Jacobian matrix of $S^j$ at $x.$

 As usual, the law of $(X_t)_{t \geq 0}$ when $X_0 = x$ is denoted $\mathbb{P}_x.$ For $x,y \in \mathbb R^n$ and $A \subset \mathbb R^n,$ we write $\langle x,y \rangle = \sum_{i = 1}^n x_i y_i, \|x\| = \sqrt{\langle x, x \rangle}$ and $d(x,A) = \inf_{y \in A} \|x-y\|.$
 If $A$ is  Borel,  we let $\tau_A = \inf \{t \geq 0 \: :X_t \in A\}$ and  $\tau_A^{out} = \tau_{\mathbb R^n \setminus A}.$

  Associated to (\ref{eq:sde}) is the Stroock and Varadhan deterministic control system: \beq
\label{controlsde}
\dot{y}(t) = S^0(y(t)) + \sum_{j = 1}^m u^j(t) S^j(y(t))
 \eeq where the control function   $u = (u^1, \ldots, u^m) : \Rp \mapsto \RR^m, $ can be chosen to be piecewise continuous.
Given such a control function, we let $y(u,x, \cdot)$ denote the  maximal solution to (\ref{controlsde}) starting from $x$ (i.e.~such
that $y(u,x,0) = x$).
  \begin{definition}
  \label{def:SDEaccessible}
 An  open set $U \subset \RR^n$  is said to be accessible  by $\{S^0, (S^j)\}$ from $x \in \RR^n$ if there exist a
   piecewise continuous control $u$ and $t \geq 0$ such that $y(u,x,t) \in U.$ If in addition, $y(u,x,s) \in \cD$ for all $0 \leq s \leq t,$ (in which case $x \in \cD$ and $U \cap \cD \neq \emptyset$) we say that $U$ is  $\cD$-accessible  by $\{S^0, (S^j)\}$ from $x.$

 By abuse of language, we say that a point $y$  is  accessible (respectively $\cD$-accessible)   by $\{S^0, (S^j)\}$ from $x,$ provided  every neighborhood $U$ of $y$ is accessible (respectively $\cD$-accessible)   by $\{S^0, (S^j)\}$ from $x.$
  \end{definition}
Throughout, we will always assume that the following hypothesis  is satisfied (or, like in Corollary \ref{cor:sde}, implied by other assumptions).
\begin{hypothesis}[H1]
The set  $\mathbb R^n \setminus \bD$ is accessible by $\{S^0, (S^j)\}$ from all $x \in \bD.$
\end{hypothesis}
It easily follows from Stroock and Varadhan support theorem and compactness of $\bD$ (see Propositions \ref{prop:control} and \ref{lem:out}) that, under  (H1),
$$\mathbb{P}_x(\tau^{out}_\bD < \infty) = 1$$ for all $x \in \bD.$

The following definition  is classical (see e.g~\cite{bass}).
\begin{definition}
\label{def:regular}
A  point $p \in \partial \cD$ is called  regular for $\mathbb R^n \setminus \bD$ with respect to  (\ref{eq:sde}), provided $\mathbb{P}_p(\tau^{out}_\bD = 0) = 1.$
\end{definition}
Our second standing hypothesis is the regularity of points in $\partial \cD.$
\begin{hypothesis}[H2]
Every point $p \in \partial \cD$ is regular for  $\mathbb R^n \setminus \bD$ with respect to  (\ref{eq:sde}).
\end{hypothesis}
For the usual Brownian motion in $\RR^n$ (i.e~$S^0 = 0, m = n,$ and $(S^1, \ldots, S^n)$ is an orthonormal basis in $\RR^n$), a
classical condition on $\cD$ ensuring (H2) is the so-called {\em exterior cone condition}: For all $p \in \partial \cD,$ there exists
an open truncated cone at vertex $p$ contained in $\mathbb R^n \setminus \bD$ (see e.g~\cite{bass}, Proposition 21.8). If the
diffusion is elliptic at $p$ (meaning that $S^1(p), \ldots S^m(p)$ span $\RR^n$), the exterior cone condition also ensures that $p$
is regular. In case the diffusion is in divergence form, this follows directly from Aronson’s
estimates~\cite{Stroock1988,LiZhang2005} which allow to bound from below the law of $(X_t)$ by (up to a constant) the law of a
Brownian motion. The general case follows by a straightforward
application of Girsanov's theorem.

For degenerate diffusions, a practical  condition   is the stronger hypothesis (H2') given below.
\begin{definition}
\label{def:unitoutward} We say that $\cD$ satisfies the exterior sphere condition at $p \in \partial \cD,$ if there exists
a vector $v \in \RR^n,$  called a  unit outward normal vector  at $p,$ such that $\|v\| = 1$ and
 $d(p + r v,\bD) = r$ for some $r > 0.$
 \end{definition}
 In other words, the open  Euclidean ball with center $p + rv$ and radius $r$ is contained in  $\mathbb R^n \setminus \bD$ and its closure touches $\bD$ at  $p.$
If the condition holds at every $p \in \partial \cD,$ we simply say that $\cD$ satisfies  the  {\em exterior sphere condition}.
 \brem
 {\rm  Examples of sets verifying the exterior sphere condition are sets with $C^2$ boundary and convex sets. In the first case every $p \in \partial \cD$ has a unique unit outward normal vector. In the second case, a point $p \in \partial \cD$ may have infinitely many  unit outward normal vectors (think of a convex polygon).
 If $\cD \subset \RR^2$ is a nonconvex polygon and $p \in \partial \cD$ is a vertex at which the interior angle is $> \pi$, the exterior sphere condition is not satisfied at $p$ although the exterior cone is.
   }
  \erem
The following condition (H2') implies (H2). This will be proved in Section \ref{sec:sde}, Proposition \ref{prop:H2'impliesH2}.
\begin{hypothesis}[H2'] $\cD$ satisfies the exterior sphere condition and for each  $p \in \partial \cD,$
 there exist an outward unit normal vector  $v$ at $p,$  and $i \in \{1, \ldots, m\}$  such that $\langle S^i(p), v \rangle \neq 0.$
\end{hypothesis}
The next result is an existence theorem.  Its proof, in
    Sections~\ref{sec:accessibility} and~\ref{sec:pf-th1.5}, shows that it only requires the Lipschitz continuity of  $S^0$ and the $C^2$ regularity of  the $S^j, j \geq 1.$ \bla{Note that this result does not require ellipticity nor hypoellipticity condition.}

\bthm
\label{th:mainresultforsde}
Assume  that:
\bdes
\iti
Hypotheses (H1) and (H2)  hold true;
  \itii For some $\eps > 0,$ the set $\cD_{\eps} = \{x \in \cD : \: d(x,\partial \cD) > \eps\}$ is $\cD$-accessible  by $\{S^0, (S^j)\}$ from all $x \in \cD \setminus \cD_{\eps}.$
  \edes
Then,   there exists   a  QSD for $(X_t)$ on $\cD.$
\ethm

Simple criteria ensuring (H1) and condition $(ii)$ of Theorem \ref{th:mainresultforsde} will be discussed in Section \ref{sec:sde}
(Propositions \ref{cor:pinsky} and \ref{cor:inward}). In
particular, Proposition \ref{cor:inward} has the following useful consequence that, when $\partial \cD$ is $C^2,$ (H2') implies
condition $(ii)$ of Theorem \ref{th:mainresultforsde}. Therefore,

\bcor Assume that Hypothesis (H1) holds true,  $\partial \cD$ is $C^2,$ and  for each $p \in \partial \cD$ one of the vectors $S^j(p), j = 1, \ldots m,$ is transverse to $\partial \cD$ (i.e~ $S^j(p) \not \in T_p \partial \cD$). Then, there exists   a  QSD for $(X_t)$ on $\cD.$
\ecor
If the SDE (\ref{eq:sde}) enjoys certain {\em hypoellipticity} properties, \bla{related to classical H\"ormander conditions as defined below}, more can be said. Given a family $\mathcal{S}$ of smooth
vector fields on $\mathbb R^n$ and $k \in \NN,$ we let $[\mathcal{S}]_k$ denote the set of vector fields recursively defined by
$[\mathcal{S}]_0 = \mathcal{S},$ and $$[\mathcal{S}]_{k+1} = [\mathcal{S}]_k \cup \{ [Y,Z] \: : Y, Z \in [\mathcal{S}]_k\}$$ where
$[Y,Z]$ stands for the Lie bracket of $Y$ and $Z.$ Set $[\mathcal{S}] = \cup_k [\mathcal{S}]_k$ and
$[\mathcal{S}](x) = \{Y(x) \: : Y \in [\mathcal{S}] \}.$
\begin{definition}
\label{def:hormanderconditions}
A  point $x^* \in \mathbb R^n$   is said to  satisfy the   weak H{\"o}rmander condition (respectively the  H\"{o}rmander condition, respectively the  strong H\"ormander condition)  if $[\{S^0, \ldots, S^m\}](x^*)$ (respectively
 $${\{S^1(x^*),  \ldots, S^m(x^*)\} \cup \{[Y,Z](x^*) :\: Y, Z \in [\{S^0, \ldots, S^m\}] \}},$$ respectively $$[\{S^1, \ldots, S^m\}](x^*)\ \text{)}$$ spans  $\RR^n.$
 \end{definition}

\bthm
\label{th:sdehypo}
Assume that:
 \bdes
\iti Hypotheses (H1) and  (H2') hold;
\itii  The weak H\"{o}rmander  condition is satisfied at every point $x \in \bD;$
\itiii Every $y \in \cD$ is  $\cD$-accessible by $\{S^0, (S^j)\}$ from every $x \in \cD.$
\edes
Then the $QSD$ $\mu$ is unique. Its  topological support equals $\bD$  and it has a smooth density with respect to the Lebesgue measure.

Suppose furthermore  that there exists a point $x^* \in \cD$ at which the weak H\"{o}rmander condition is strengthened to  the
H\"{o}rmander condition. Then there exist $\alpha > 0$, $C \in ]0,+\infty[$ and a continuous function $h: \cD \mapsto ]0,\infty[$ with
$h(x) \rar 0,$ as $x \rar \partial \cD,$ satisfying $\mu(h) = 1$, such that for all $\rho \in \cM_1(\cD)$ (the set of probability measures over $\cD$),
\begin{equation}
  \label{eq:cv-sde}
  \|\mathbb{P}_{\rho}(X_t \in \cdot \mid \tau_\cD^{out} > t) - \mu(\cdot)\|_{TV} \leq \frac{C}{\rho(h)} e^{-\alpha t}
\end{equation}
where $\| \cdot \|_{TV}$ stands for the total variation distance and $\rho(h) := \int h d\rho.$
\ethm

Note that Condition (ii) of Theorem~\ref{th:mainresultforsde} readily follows from Condition (iii) of
  Theorem~\ref{th:sdehypo}, so that the assumptions of Theorem~\ref{th:sdehypo} imply the existence of a QSD.

\bcor \label{cor:sde} Assume that Hypothesis (H2') holds, and  that the strong H\"{o}rmander  condition is satisfied at every point
$x \in \bD.$ Then, the conclusions of Theorem \ref{th:sdehypo} hold true and the density of  $\mu$ is positive on $\cD.$ \ecor
\prf Let, for $\eps \geq 0,$ $y^{\eps}(x,u,\cdot)$ be defined like $y(x,u,\cdot)$ when $S^0$ is replaced by $\eps S^0.$ By Chow's theorem  (see e.g~ \cite{Jurdjevic}, Chapter 2, Theorem 3), the  strong H\"{o}rmander  condition implies that for all $x,y \in \cD$ there exists a control $u$ piecewise continuous with $u^i(s) \in \{-1,0,1\}$ and $t \geq 0$ such that $y^0(x,u,s) \in \cD$ for all $0 \leq s \leq t$ and $y^0(x,u,t) = y.$  By continuity of solutions to an ODE, with respect to some parameter (or by a simple application of Gronwall's lemma), $y^{\eps}(x,u,\cdot) \rar y^0(x,u,\cdot)$ uniformly on $[0,t]$ as $\eps \rar 0.$ Let $u^{\eps}(s) = \frac{u(\frac{s}{\eps})}{\eps}.$  Then
$y^{\eps}(x,u,s) = y(x,u^{\eps},\eps s).$ This proves that $y$ is $\cD$-accessible from $x.$ Because the strong H\"{o}rmander condition
also holds in a neighborhood of $\bD,$ the same proof also shows that, for $y$ in a neighborhood of $\bD$ and $x \in \bD,$ $y$ is
accessible from $x$, so (H1) is satisfied. The conditions, hence the conclusions, of Theorem \ref{th:sdehypo} are then satisfied. Positivity of $\frac{d\mu}{dx}$ is proved in Lemma \ref{lem:hyposmooth}. \qed
\brem
\label{rem:manifold} {\rm  The results above easily extend to the situation where $\mathbb R^n$ is replaced by a $n$-dimensional manifold, provided $\partial \cD$ is a $(n-1)$-dimensional $C^2$ sub-manifold. This is illustrated in Examples \ref{ex:cylinder} and \ref{ex:cylinder2}.}
\erem
\brem {\rm The  paper \cite{GNW2020} considers QSDs and their properties under the assumption that the underlying process is
  strong Feller. Quasi-stationarity for strong Feller diffusion processes have also been a focus of interest in the
  recent literature \cite{GNW2020, LRR2021, R2021}. Observe that none of the conditions of Theorem \ref{th:mainresultforsde}, neither
  the weak H\"{o}rmander condition assumed in Theorem \ref{th:sdehypo} imply that the semigroup induced by  (\ref{eq:sde}) is strong
  Feller (see Example~\ref{ex:cylinder} below).} \erem
\bex
\label{ex:cylinder}
{\rm
 We consider the situation where the SDE is defined on the cylinder $\RR/ \ZZ \times \RR$ (instead of $\mathbb R^n$). Let $m = 1,$  $$S^0(x,y) = \partial_x \mbox{ and } S^1(x,y) = \partial_y.$$ Let $\cD = \RR/ \ZZ \times ]0,1[.$ Here the conditions of Theorem \ref{th:sdehypo}, are easily seen to be satisfied. However, the process is not strong Feller (the dynamics in the $x$-variable being deterministic). It is not hard to check that the unique QSD is the measure
\beq
\label{eq:QSDcylinder}
\mu(dx dy) = 2\frac{\sin(\pi y)}{\pi} \Ind_\cD(x,y) dx dy.
\eeq
} \eex
\bex[Example \ref{ex:cylinder}, continued]
\label{ex:cylinder2}
{\rm Consider the setting of Example \ref{ex:cylinder}, but with $$S^0(x,y) = a(y)\partial_x \mbox{ and } S^1(x,y) = \partial_y,$$ where $a$ is a smooth function $\geq 1.$ Like in example \ref{ex:cylinder}, the unique QSD $\mu$ is given by (\ref{eq:QSDcylinder}). Suppose that  $a'(y^*) \neq 0$ for some  $0 < y^* < 1.$ Then  the H\"{o}rmander condition holds at $(x,y^*)$, so that by Theorem \ref{th:sdehypo},  $(\mathbb{P}_{x}(X_t \in \cdot | \tau_\cD^{out} > t))_{t \geq 0}$ converges at an exponential rate to $\mu.$}
\eex
\brem {\rm In absence of condition $(iii)$ in Theorem \ref{th:sdehypo}, there is no guarantee that the QSD is unique, as shown by the
  next example. Still every QSD has a smooth density (see Lemma \ref{lem:hyposmooth}).} \erem
\bex {\rm
    Let $n=1$ and $\cD=]0,5[$. We consider smooth functions $\varphi^1,\varphi^2,\psi^1,\psi^2:\cD\to[0,1]$ such that $\psi^1+\psi^2=1$ and
    \begin{align*}
    \begin{cases}
    \varphi^1_{\rvert ]0,1]}\equiv 1,\ \varphi^1_{\rvert [1,2]}\leq 1,\ \varphi^1_{\rvert [2,5[}\equiv 0,\\
    \varphi^2_{\rvert ]0,3]}\equiv 0,\ \varphi^2_{\rvert [3,4]}\leq 1,\ \varphi^2_{\rvert [4,5[}\equiv 1,\\
    \psi^1_{\rvert ]0,2]}\equiv 1,\ 0<\psi^1_{\rvert  ]2,3[}\leq 1,\ \psi^1_{\rvert [3,5[}\equiv 0,\\
    \psi^2_{\rvert ]0,2]}\equiv 0,\ 0<\psi^2_{\rvert  ]2,3[}\leq 1,\ \psi^2_{\rvert [3,5[}\equiv 1.
    \end{cases}
    \end{align*}
    For all $\alpha>0$, we consider the absorbed diffusion process $X^\alpha$ evolving according to the It\^o SDE
    \begin{align*}
    dX^\alpha_t=\left(\varphi^1(X^\alpha_t)+\sqrt{\alpha}\,\varphi^2(X^\alpha_t)\right)\, dB_t\,+(\psi^1(X_t)+\alpha\,\psi^2(X_t))\,dt.
    \end{align*}
    and absorbed when it reaches $\partial \cD=\{0,5\}$. While $X^\alpha$ satisfies the conditions (i) and (ii) of Theorem \ref{th:sdehypo},  it does not satisfy condition (iii) and it admits either one or two QSDs, depending on the value of $\alpha$. Indeed, as shown in~\cite{BenaimChampagnatEtAl}, there exists $\alpha_c>0$ such that
    \begin{itemize}
        \item for all $\alpha\in]0,\alpha_c]$, $X^\alpha$ admits a unique QSD supported by $[3,5]$,
        \item for all $\alpha\in]\alpha_c,+\infty[$, $X^\alpha$ admits exactly two QSDs, supported respectively by $[3,5]$ and $[0,5]$.
    \end{itemize}
} \eex
\section{Killed processes}
\label{sec:main}
Throughout, we let $M$ denote a separable and locally compact metric  space,  and $D \subset M$  a nonempty set with compact closure
 $K = \overline{D}$ in $M$, such that $D$ is  open {\em relative to} $K.$ That is
 $D = {\cal O} \cap K$ for some open set ${\cal O} \subset M.$  Considering such a general setting is important for
   example for piecewise deterministic Markov processes (see Section~\ref{sec:PDMP}; see also Example 2.6 below).


\bla{\subsection{Preliminary results}
All the definitions and preliminary results gathered in this section are classical, or adapted from classical ones.}

 We let $(P_t)_{t \geq 0}$ denote a Markov Feller semigroup\footnote{\bla{We also mentioned strong Feller semigroups in the
     introduction. They are defined as Feller semigroups, except for the property that $P_t$ maps $C_0(M)$ to itself, which is
     replaced by the property that $P_t$ maps bounded measurable functions on $M$ to bounded continuous functions on $M$.}} on $M$. By
 this, we mean (as usual) that $(P_t)_{t \geq 0}$ is a semigroup of Markov operators on $C_0(M)$ (\bla{the closure for the supremum
   norm of the set of continuous, compactly supported functions from $M$ to $\mathbb{R}$}) and that $P_t f(x) \rar f(x)$ as $t \rar 0$ for the uniform norm for all $f \in C_0(M)$. Observe that
 since we are interested by the behavior of the process killed outside $K$, the behavior of $(P_t)_{t \geq 0}$ at infinity is
 irrelevant (and the reader can think of $M$ as compact without loss of generality). Feller processes, named after William Feller,
 forms\ a large class of processes whose sample path properties have been studied, among others, by Doob, Kinney, Dynkin, Blumenthal,
 Hunt and Ray.
   We refer the reader to~\cite{Chung1982,Kallenberg2002} for properties and historical notes on these processes, see also~\cite{Legall2,KaratzasShreve1991,RogersWilliams2000} for expository texts on Feller semi-groups and their sample path properties.

 By classical results (see e.g~Le Gall \cite{Legall2}, Theorem 6.15),  there exist a filtered space  $(\Omega, {\cal F}, ({\cal F}_t))$ with  $({\cal F}_t)$  right continuous and complete, a family of probabilities  $(\mathbb{P}_x)_{ x \in M}$ on $(\Omega, {\cal F})$ and a  continuous time adapted process $(X_t)$ on $(\Omega, {\cal F}, ({\cal F}_t))$ taking values in $M,$  such  that:
\bdes
\iti $(X_t)$ has cad-lag paths,
\itii $\mathbb{P}_x(X_0 = x) = 1$ and,
 \itiii $(X_t)$ is a Markov process with semigroup $(P_t),$ meaning that
 $$\mathbb{E}_x(f(X_{t+s})|{\cal F}_t) = P_s f(X_t)$$ for all $t,s \geq 0$ and $f$ measurable bounded (or $\geq 0$).
 \edes
For any Borel set $A \subset M$ we let $\tau_A = \inf \{t \geq 0: \: X_t \in A\}$ and $\tau^{out}_A = \tau_{M \setminus A}.$
The assumptions on $({\cal F}_t)$ (right continuous and complete) imply that $\tau_A$ and $\tau^{out}_A$ are stopping times with respect to  $({\cal F}_t)$
(see e.g.~Bass \cite{10.1214/ECP.v15-1535}).
\brem
 {\rm  For $X_0 = x \in D,$ $\tau^{out}_D \leq \tau^{out}_{K}$ but it is not true in general  that  $\tau^{out}_D = \tau^{out}_{K}.$ Consider for example the ODE on $\RR^2$ given by
                  $$ \left\{ \begin{array}{c}
                      \dot{x} = 1 \\
                       \dot{y} = 0
                   \end{array} \right. $$
Let $$D = \{(x,y) \in \RR^2 \: -1 < x< 1, - 1 < y < x^2 \}.$$
 For $-1 < x < 0$ and $y = 0$ the trajectory $(x(t),y(t)) = (x +t,0)$ starting at $(x,0)$ leaves $D$ at time $-x$ and $K$ at time $-x +1.$ }
\erem
 We now recall the definition of accessibility for general Feller processes.
\begin{definition}
\label{def:accessible}
An open set $U \subset M$ is said to be accessible (by $(P_t)$) from $x \in M,$  if there exists $t \geq 0$ such that $P_t(x,U) = P_t \Ind_U(x) > 0.$

Here also, by abuse of language, we say that  $y \in M$ is   accessible from $x \in M,$  if every open neighborhood $U$ of $y$ is accessible from $x.$
\end{definition}

We shall assume throughout all this section that the following assumption is satisfied.
\begin{hypothesis}[Standing Hypothesis]
\label{hyp:main} The set $M \setminus K$ is accessible  from all $x \in K.$ \end{hypothesis}

The next result is a basic, but very useful, consequence of this  assumption.
It is interesting to point out that it only requires  Feller continuity, accessibility of $M \setminus K$ and compactness of $K$.
 Note that similar - albeit different- results  can be found in the literature (see e.g~\cite{bass}, Proposition 21.2).

\bprop
\label{lem:out}
  There exist positive constants $C, \Lambda$ such that $$\mathbb{P}_x(\tau^{out}_K > t) \leq C e^{-\Lambda t}$$ for all $t\geq 0$ and $x \in K.$ In particular
 $\tau^{out}_K < \infty$, and hence $\tau^{out}_D < \infty$, $\mathbb{P}_x$ almost surely for all $x \in K.$
\eprop
\prf Let $U=M\setminus K$. By Feller continuity, for all $t \geq 0,$ the set $O_t = \{x \in M: \: P_t(x,U) > 0\}$ is open (possibly empty).
We detail the proof since the argument will be used several times in the sequel. We notice that the sequence of
  functions $f_n(x)=(nd(x,M\setminus U)\wedge 1)$ is non-decreasing and converges to $\mathbbm{1}_U(x)$. Therefore, $x\in
  O_t$ if and only if $x\in O_t^n= \{x \in M: \: P_t f_n(x) > 0\}$ for some $n\geq 1$. Since $O^n_t$ is open by Feller continuity
  of the semi-group $P_t$, we deduce that $O_t$ is open. By  the standing hypothesis, the family $\{O_t : \: t \in \RR^+\}$ covers $K,$ so that,
by compactness, there exist $t_1, \ldots, t_n$ such that $K \subset \cup_{i = 1}^n O_{t_i}.$ In particular, for some $\delta > 0$ and
$t = \max \{t_1, \ldots, t_n \}$ $\mathbb{P}_x(\tau_U > t) \leq 1-\sup_{1\leq i \leq n}P_{t_i}(x,U)\leq 1-\delta$ for all $x \in K.$ Thus, by the Markov property,  $\mathbb{P}_x(\tau_U > k t) \leq (1-\delta)^k.$ This proves the result. \qed

 \subsection{Green kernel and QSDs}
Let $B(D)$  denote the set of bounded  measurable functions $f : D \mapsto \RR.$ For all $f \in B(D), x \in D$  and $t \geq 0$ set
 $$P_t^D f(x) = \mathbb{E}_x(f(X_t) \Ind_{\tau^{out}_D > t}).$$
Then, $(P_t^D)_{t \geq 0}$
is a well defined  sub-Markovian semigroup on $B(D).$  The semigroup property is a consequence of the Markov property and the fact
that $\tau^{out}_D$ is a first exit time and a stopping time~\cite[Section 2.6]{collet2012quasi}.

We define as usual the {\em Green kernel} $G^D$ as the bounded (by
Proposition~\ref{lem:out}) operator defined on $B(D)$ given by
$$G^D f(x) = \int_0^{\infty} P_t^D f(x) dt = \mathbb{E}_x \left(\int_0^{\tau^{out}_D} f(X_t) dt\right).$$
For all $x \in D$ and $A \subset M,$ a Borel set, we let $$G^D(x,A) = G^D \Ind_{A \cap D} (x).$$
\subsubsection*{Quasi-stationary distributions}
A {\em quasi-stationary distribution} (QSD) for $(P_t^D)$ is a probability measure $\mu$ on $D$ such that, for all $t\geq 0$,
\beq
\label{eq:QSDdef}
\mu P_t^D = e^{-\lambda t} \mu
\eeq for some $\lambda > 0.$ For further reference we call $\lambda$ the {\em absorption rate} (or simply the rate) of $\mu.$
Equivalently~\cite{Meleard-Villemonais-2012},
$$\frac{\mu P_t^D(\cdot)}{\mu P_t^D \Ind_D} = \mathbb{P}_{\mu}(X_t \in \cdot \mid \tau_D^{out} > t) = \mu(\cdot).$$

\blem
\label{lem:QDS}
Equation (\ref{eq:QSDdef}) holds if and only if $\mu G^D = \frac{1}{\lambda} \mu.$
\elem
\prf
Clearly, by definition of $G^D$, equation (\ref{eq:QSDdef}) implies that $\mu G^D = \frac{1}{\lambda} \mu.$ Conversely, assume that  $\mu G^D = \frac{1}{\lambda} \mu.$ Then for every bounded nonnegative measurable map $f : D \mapsto \RR,$ $\mu(G^D f) = \frac{1}{\lambda} \mu(f)$ and also
$$\mu (G^D P_t^D f) =  \frac{1}{\lambda} \mu(P_t^D f).$$ That is $$\mu(\int_t^{\infty} P_s^D f ds) = \frac{1}{\lambda} \mu(P^D_t f).$$ Equivalently,
$$\mu (G^D f - \int_0^t P_s^D f ds) = \frac{1}{\lambda} \mu (P_t^D f).$$ This shows that the bounded map $v(t) = \mu(P_t^D f)$ satisfies the integral equation $$v(t) - v(0) = - \lambda \int_0^t v(s) ds,$$ for all $t\geq 0$. It follows that $v(t) = v(0) e^{- \lambda t}.$
\qed
Let $C_b(D) \subset B(D)$ denote the set of bounded continuous functions on $D,$ and $C_0(D) \subset C_b(D)$ {\color{black} the
  closure of the set of continuous, compactly supported functions from $D$ to $\mathbb{R}$. Observe that, since $K$ is compact,
  $C_0(D)$ is also the subset of functions $f$ such that $f(x) \rar 0$ when $d(x,\partial_K D)\to 0$, where $d(x,\partial_K D)$ is
  the distance between $x$ and $\partial_K D:=K\setminus D$. 
  Similarly, we define $C_b(K)$ as the set of bounded continuous functions on
  $K$.}

\brem{\rm In the recent paper~\cite{FerreRoussetStoltz2021}, the authors prove the existence and convergence to a QSD under the condition that the sub-Markovian semigroup is strong Feller.
Observe that in our case, although $(P_t)$ is Feller, there is no evidence in general that  $(P_t^D)$ is strong Feller nor that it preserves $C_b(D).$ On the other hand, under rather weak, reasonable conditions, $G^D$ maps $C_b(D)$ into $C_0(D)$, as illustrated by the following example.

}\erem
 \bex
 {\rm Consider the ODE on $\RR$ given by $\dot{x} = -1.$  For $D = ]0,1[$ $P_t^Df(x) = f(x-t) \Ind_{x > t}$ is not Feller,
 but $G^D f(x) = \int_0^x f(u) du$ is Feller (and even strong Feller). If now $D = ]0,1],$ then \bla{$D$ is open relative to
   $M=]-\infty,1]$, with compact closure and} $G^D$ maps $C_b(D)$ into $C_0(D).$
}
 \eex
The condition that $G^D(C_b(D)) \subset C_0(D)$ plays a key role in the next  theorems and will be
    investigated in the subsequent sections for degenerate diffusions and PDMPs.

\begin{definition}
\label{def:Daccessible}
 An open set $U \subset M$ is said  $D$-accessible (by $(P_t)$) from $x \in D$ if $P_t^D(x,U) > 0$ for some $t \geq 0.$
 A  point $y \in K$ is said  $D$-accessible from $x \in D$ if every open neighborhood of $y$ is $D$-accessible from $x.$
\end{definition}
\blem
\label{lem:GDaccess}
An open set $U \subset M$ is $D$-accessible from $x \in D$ if and only if $G^D(x,U) > 0.$ In particular, the set of $D$-accessible points from $x$ coincides with the topological support of $G^D(x,\cdot).$
\elem
\prf Suppose that $P^D_t(x,U) > 0$ for some open set $U$ and $t \geq 0.$   Fatou Lemma and right continuity of paths imply that
$$\liminf_{s \rar t, s > t} P^D_s(x,U) \geq \mathbb{E}_x(\liminf_{s \rar t, s > t} \Ind_U(X_s) \Ind_{\tau_D^{out} > s}) \geq P^D_t(x,U) > 0.$$
This proves that $s \rar P^D_s(x,U)$ is positive on some interval $[t,t+\eps],$ hence $G^D(x,U) > 0.$ The converse implication is obvious. \qed

We now state and prove our first main result.
\bthm
\label{th:main}
Assume that $G^D(C_b(D)) \subset C_0(D).$ Then, the  conditions (i),(ii) below are equivalent and imply the existence of a QSD.
\bdes
\iti There exists an open set $U \subset M,$  $D$-accessible from all $x \in D$ and such that $ \overline{U \cap K} \subset D.$
\itii  There exists a function $\phi \in C_0(D),$  positive on $D$ and $\theta > 0,$ such that
$G^D \phi \geq \theta \phi.$
\edes
\ethm

\medskip\noindent\prf
We first show that  $(i)$ implies $(ii).$ For $\eps > 0,$ let $$\overline{U}^{\eps}_K = \{x \in K : \: d(x,U \cap K) \leq \eps\}.$$ Choose $\eps > 0$ small
enough so that $\overline{U}_K^{\eps}  \subset D$.

Let $\psi(x) = (1- \frac{d(x, U \cap K)}{\eps})^+$ and $\phi = G^D \psi.$ Then, $\phi\in C_0(D)$ (because
$G^D(C_b(D))\subset C_0(D)$) and for all $x \in D,$ $$\phi(x) \geq G^D \Ind_{U \cap K}(x) = G^D \Ind_{U}(x) > 0,$$ where the last inequality holds by Lemma~\ref{lem:GDaccess}. Thus, by compactness of
$\overline{U}^{\eps}_K,$ $\theta = \inf_{x \in \overline{U}_K^{\eps}} \phi(x) > 0.$ Since $\psi = 0$ on $K \setminus
\overline{U}_K^{\eps},$ it follows that for all $x \in D,$ $$\phi(x) \geq \theta
  \mathbbm{1}_{\overline{U}^\varepsilon_K}\geq \theta \psi(x).$$ Therefore,  $G^D \phi(x) \geq \theta \phi(x).$

We now prove that $(ii)$ implies $(i)$. For all $\eps \geq 0$ and $ x_0 \in D,$ set  $$V_{\eps} = \{y \in D:\: \phi(y) > \eps\},$$ and
 $$\eta = \eta(x_0) := \inf \{\eps \geq 0 : \: G^D(x_0,V_{\eps}) = 0\}.$$ Note that
   $V_{\|\phi\|}=\emptyset$, so that $\eta<+\infty$. By monotone convergence, $G^D(x_0,V_{\eta}) = \lim_{n \rar \infty} G^D(x_0,V_{\eps+1/n}) = 0.$  In particular, $\eta$ is positive because, by (ii), $G^D(x_0,D)$ is. Therefore, by definition of $\eta,$ $G^D(x_0,V_{\eta/2}) > 0,$ so that
$\mathbb{P}_{x_0} ( X_t \in V_{\eta/2},\,t<\tau_D^{out}) > 0$ for some $t \geq 0.$ Furthermore,  $G^D(X_t, V_{\eta}) = 0$ $\mathbb{P}_{x_0}$-almost surely on the event $\{t < \tau_D^{out} \}.$ Indeed,
\begin{multline*}
  \mathbb{E}_{x_0} (G^D(X_t,V_{\eta}) \Ind_{t < \tau_D^{out}}) = P_t^D G^D(\Ind_{V_{\eta}})(x_0) \\ = \int_{t}^{\infty} P_s^D(x_0,V_{\eta})ds
  \leq G^D(x_0,V_{\eta}) = 0.
\end{multline*}
In particular, $X_t \in D \setminus V_{\eta}$ for all $\tau_D^{out} > t \geq 0$ almost surely, by Lemma \ref{lem:GDaccess}.
This shows that there exists a point $y \in V_{\eta/2} \setminus V_{\eta}$ with $G^D(y,V_{\eta}) = 0.$  Thus
\begin{equation}
  \label{eq:pf-Thm2.9}
  \theta \leq \frac{G^D \phi (y)}{\phi(y)} = \frac{G^D (\phi \Ind_{D \setminus V_{\eta}}) (y)}{\phi(y)}
  \leq 2 G^D \Ind_D(y) \leq 2\sup_{x \in D \setminus V_{\eta}  }G^D \Ind_D(x).
\end{equation}
Now, since $G^D\Ind_D\in C_0(D)$, there exists $\alpha>0$ such that
\[
\sup_{x\in D,\, d(x, \partial_K D)\leq\alpha} G^D\Ind_D(x)<\frac{\theta}{2}
\]
and, since $\phi$ is positive on $D$ and $\phi\in C_0(D)$, there exists $\eta_0>0$ such that $D\setminus V_{\eta_0}\subset\{x\in D,\,
d(x, \partial_K D)\leq\alpha\}$. Therefore,~\eqref{eq:pf-Thm2.9} implies that  $\eta(y)>\eta_0$, hence
  $\eta=\eta(x_0)>\eta_0$. Since $x_0$ was arbitrary we deduce that the map $x\mapsto\eta(x)$ is bounded from below by
$\eta_0$. Since $G^D({\color{blue}x_0}, V_{\eta_0/2})\geq G^D(x_0,V_{\eta/2})>0$, Lemma~\ref{lem:GDaccess} proves
$(i)$ with $U = V_{\eta_0/2}$ (and actually $\overline{U}\subset D$).

Our last goal is to prove the existence of a QSD under Condition(ii).
Let $B(K)$ (respectively $C_b(K)$) be the space of bounded (respectively   bounded continuous) functions over $K.$ The operator $G^D$ extends to a bounded operator  $G^{D, K}$ on $B(K)$ defined as
$$ {G^{D,K}}f(x) = \left \{ \begin{array}{ll}
  G^D (f|_D)(x) & \mbox{ for }  x \in D \\
  0  & \mbox{ for }  x \in \partial_K D
\end{array} \right.$$
The assumption that $G^D(C_b(D)) \subset C_0(D)$ implies that   ${G^{D,K}}(C_b(K)) \subset C_b(K).$

Let ${\cal M}_1(\phi)$ be the set of Borel finite nonnegative measures $\mu$ on $K$ such that $\mu(\phi) = 1$ (extending
  $\phi$ to $K$ by setting $\phi(x)=0$ for all $x\in\partial_K D$) and let $T : {\cal M}_1(\phi) \mapsto {\cal M}_1(\phi)$ be the map defined by
$$T(\mu) = \frac{\mu {G^{D,K}}}{\mu {G^{D,K}} \phi}$$
We first observe that $T$ is continuous for the weak* topology: if $\mu_n\rightharpoonup \mu$ for the weak* topology in
$\mathcal{M}_1(\phi)$, then $\mu_n {G^{D,K}}\rightharpoonup \mu {G^{D,K}}$ and $\mu_n {G^{D,K}}\phi\rightarrow\mu {G^{D,K}}\phi$. Since in addition $\mu_n
{G^{D,K}}\phi\geq\theta\mu_n(\phi)=\theta>0$ by $(ii),$ $T(\mu_n)\rightharpoonup T(\mu)$.

Choose an open  neighborhood ${\cal N}$ of $\partial_K D$ such that ${G^{D,K}} \Ind_K \leq \theta/2$ on ${\cal N} \cap K$ and set $C = \sup_{x \in K \setminus {\cal N}} \frac{({G^{D,K}} \Ind_K)(x)}{\phi(x)} < \infty.$ Then, for all $\mu \in {\cal M}_1(\phi)$
$$T(\mu)(K) = \frac{\mu ({G^{D,K}} \Ind_K)}{\mu {G^{D,K}} \phi} \leq \frac{\mu ({G^{D,K}} \Ind_K)}{\theta} \leq \frac{ \theta/2 \mu(K) + C \mu(\phi)}{\theta} = \frac{\mu(K)}{2} + \frac{C}{\theta}.$$
It follows that for any $R \geq 2 \frac{C}{\theta}$ the  set $${\cal M}_1^R(\phi) = \{\mu \in {\cal M}_1(\phi) \: : \mu(K) \leq R\}$$ is invariant
 by $T.$ Since ${\cal M}_1^R(\phi)$ is convex and compact (for the weak* topology),  $T$ admits a fixed point by Tychonov's Theorem.
 If $\mu$ is such a fixed point, $\mu(D)\leq R$ and $\mu(D)>0$ since $\mu(\phi)=1$ and $\phi(x)=0$ for all $x\in K\setminus D$, so we
   can define the probability measure $\bar{\mu}=\frac{\mu(\cdot \cap D)}{\mu(D)}$ on $D$. Since, for all $f\in B(K)$,
   $\mu {G^{D,K}} f=\mu(\Ind_D G^D f|_D)$, we deduce that, for all $f\in B(D)$, $\bar{\mu}G^D f=(\mu {G^{D,K}}\phi)\bar{\mu}(f)$ with $\mu {G^{D,K}}\phi>0$.
   Then $\mu$ is a QSD for $(P_t^D)$ by Lemma~\ref{lem:QDS}. \qed
 \brem {\rm The last part of the proof of Theorem \ref{th:main} is reminiscent of the proof of the existence Theorem 4.2 in \cite{ColMartMel}.} \erem
\subsection{Uniqueness  and convergence criteria}
\subsubsection*{Right eigenfunctions}
We say that  $h\in B(D)$ is a {\em positive right eigenfunction} for $G^D$ if $h(x) > 0$ for all $x \in D$ and
\beq
\label{eq:righteigen}
G^D h = \frac{1}{\lambda} h
\eeq
for some $\lambda > 0.$

\blem
\label{lem:righteigen}If $h$ is a positive right eigenfunction  for $G^D$, the parameter $\lambda$ in (\ref{eq:righteigen}) necessarily equals the absorption rate of any QSD
and $$P_t^D h = e^{-\lambda t} h$$ for all $t \geq 0.$
 \elem
 \prf If $\mu$ is a QSD with rate $\lambda'$, then $\mu G^D h = \frac{1}{\lambda'} \mu h = \frac{1}{\lambda} \mu h$ so that $\lambda = \lambda'.$
 The proof of the second statement is similar to the proof of Lemma \ref{lem:QDS} and left to the reader. \qed
 Observe that there is no assumption in Lemma \ref{lem:righteigen} that $G^D(C_b(D)) \subset C_0(D)$ but when this is the case,
if there exists a positive right eigenfunction for $G^D$, then there always exists a QSD,
 by application of Theorem \ref{th:main}.

The following result shows that a strengthening in the assumptions of Theorem \ref{th:main}  ensures the existence of a positive right eigenfunction for $G^D.$
 \bcor
 \label{cor:main}
 Assume that:
\bdes
\iti $G^D(C_b(D)) \subset C_0(D)$ and $G^D$ is a compact operator \bla{from $C_0(D)$ to itself, where $C_0(D)$ is endowed with the uniform distance};
 \itii For all $x,y \in D$ $y$ is $D$-accessible from $x.$
 \edes
 Then, there exists a positive right eigenfunction $h\in C_0(D)$ for $G^D.$
 \ecor
 \prf Let $r = \lim_{n \rar \infty} \|(G^D)^n\|^{1/n}$ be the spectral radius of $G^D$ on $C_0(D).$ Let $\mu$ be a QSD (whose
 existence is given by Theorem \ref{th:main}) with rate $\lambda.$   For any $f \in C_0(D)$ such that $0 \leq f \leq 1$ and $\mu(f) \neq 0$
 $\|(G^D)^n\| \geq \mu((G^D)^n f) = \frac{1}{\lambda^n} \mu(f).$ Hence $r \geq \frac{1}{\lambda} > 0.$

 Let  $C_0^+(D) = \{f \in C_0(D) \: : f \geq 0\}.$ It is readily  seen that $C_0^+(D)$ is a {\em reproducing} cone in $C_0(D)$ invariant by $G^D,$ meaning that $C^+_0(D)$ is a cone,  $C_0(D) = \{u-v: \: u,v \in C^+_0(D)\}$ and $G^D(C^+_0(D)) \subset C^+_0(D).$ Therefore, by Krein Rutman Theorem (\cite{Deimling}, Theorem 19.2), there exists $h \in C^+_0(D) \setminus \{0\}$ such that
 $G^D h= r h.$
 Let $y \in D$ be such that $h(y) > 0.$ Then, $h \geq \frac{h(y)}{2} \Ind_{U \cap D}$ for some neighborhood $U$ of $y.$ Therefore, by
 $(ii),$ for all $x\in D,$
 $$r h(x) = (G^D h) (x) \geq  \frac{h(y)}{2}(G^D\Ind_{U \cap D})(x) > 0.$$
 This concludes the proof. \qed

\subsubsection*{Uniqueness and convergence}
Similarly as in the conservative setting, we say that $(P_t^D)$ is {\em irreducible} if there exists a nontrivial positive measure $\xi$ on $D$ such that for all $x \in D$ and $A$ Borel,
$$\xi(A) > 0 \Rightarrow G^D(x,A) > 0.$$

\bthm
\label{th:uniqueness} Assume that:
\bdes
\iti $G^D(C_b(D)) \subset C_0(D)$;
 \itii There exists a positive right eigenfunction $h\in B(D)$ for $G^D;$
 \itiii  $(P_t^D)$ is  irreducible.
  \edes
  Then $(P_t^D)$ has a unique QSD. \ethm
  \brem
  \label{rem:oncormain}
{\rm Conditions $(i)$ and $(ii)$ of Theorem \ref{th:uniqueness} are implied by the assumptions of Corollary
  \ref{cor:main}.}
\erem
\prf Let $h$ be a positive right eigenfunction and $\mu$ a QSD with rate $\lambda.$ Let $Q$ and $\pi$ respectively denote the Markov kernel and the probability  on $D$  defined by
$$Q(f) = \lambda \frac{G^D (fh)}{h},$$ and
$$\pi(f) = \frac{\mu(fh)}{\mu(h)}$$ for all $f \in B(D).$ Then, $\pi$ is invariant by $Q.$ The assumption that $(P_t^D)$ is irreducible
makes $Q$ irreducible, in the sense that $\xi(A) > 0 \Rightarrow Q(x,A) > 0$ for all $x \in D$ and $A$ Borel.
Therefore, by a standard result (see~e.g \cite{duf00} or \cite{MT93}), $\pi$ is the unique invariant probability of $Q.$ Assume now
that $\nu$ is another QSD with rate $\alpha.$ Then $\nu (G^Dh) = \frac{1}{\alpha} \nu (h) = \frac{1}{\lambda} \nu(h).$
Since $h$ is positive on $D$, we deduce that $\alpha = \lambda.$ It follows that the probability $\pi'$ defined like $\pi$ with $\nu$ in place of $\mu$ is invariant by $Q.$ By uniqueness, $\pi = \pi'$ and consequently $\mu = \nu.$ \qed
A sufficient (and often more tractable than the definition) condition  ensuring irreducibility is given by the next lemma.
\blem
\label{lem:irreduc}
Suppose that there exists an open set $U \subset M,$ $D$-accessible from all $x \in D,$ and a non trivial measure $\xi$ such that for all $x \in U$
$G^D(x,\cdot) \geq \xi(\cdot).$
Then $(P_t^D)$ is irreducible.
\elem
\prf For all $x \in D,$ there exists, by $D$-accessibility, $t \geq 0$ such that $P_t^D(x,U) > 0.$ For every Borel set $A \subset M,$
$$G^D(x,A) \geq \int_0^{\infty} P_{t+s}^D(x,A) ds = \int_0^{\infty} \int_M P_t^D(x,dy) P_s^D(y, A) ds$$
$$\geq \int_0^{\infty} \int_U P_t^D(x,dy) P_s^D(y, A) ds \geq P_t(x,U) \xi(A).$$
\qed
 If  the local minorization  $G^D(x,\cdot) \geq \xi(\cdot)$ appearing in Lemma \ref{lem:irreduc} can be improved to local minorization involving $(P_t^D),$ we also get  the exponential convergence of the conditional semigroup toward $\mu.$ More precisely,
\bthm
\label{th:convergence}
Suppose that:
\bdes
\iti $G^D(C_b(D)) \subset C_0(D)$;
 \itii There exists a positive right eigenfunction $h\in C_0(D)$ for $G^D$;
 \itiii There exist an open set $U \subset M,$ $D$-accessible from all $x \in D,$  a non trivial measure $\xi$ on $D,$ and  $T> 0$ such that for all $x \in U,$ $$P_T^D(x, \cdot) \geq \xi(\cdot).$$
 \edes Then
there exist $C, \alpha > 0$ such that, for all $\rho \in \cM_1(D)$, 
  $$\left\|\frac{\rho P_t^D}{\rho P_t^D \Ind_D} - \mu(\cdot)\right\|_{TV} \leq \frac{C}{\rho(h)} e^{-\alpha t}$$ where $\| \cdot \|_{TV}$ stands for the total variation distance.
  \ethm
\prf Let $h$ be a positive right eigenfunction and $\mu$ a QSD with rate $\lambda.$ For all $f \in B(D)$ and $t \geq 0,$ let
$$Q_t(f) = e^{ \lambda t} \frac{ P_t^D (f h)}{h}$$
and
$$\pi(f) = \frac{\mu(f h)}{\mu(h)}$$
It is readily seen that $(Q_t)$ is a Markov semigroup (usually called the $Q$--process induced by $\mu$)  having $\pi$ as invariant probability.
In order to prove the theorem we will show that:
\bdes
\item{{\bf Step 1}} There exists a probability $\nu$ on $D$ and $T_0>0$ such that for every compact set $\tilde{K} \subset D$ there is some integer $n$ and some constant $c$ (both depending on $\tilde{K}$) such
 that $$Q_{n T_0}(x,\cdot) \geq c \nu(\cdot)$$ for all $x \in \tilde{K}.$
 \item{{\bf Step 2}} The function $V = 1/h$ is a continuous and proper {\em Lyapunov function} for $Q_{T_0}$, that is
$$\lim_{x \rar  \partial_K D} V(x) = \infty,$$ and
\beq
\mylabel{eq:VLyap}
Q_{T_0} V \leq \rho V + C
\eeq for some $0  \leq \rho < 1$ and $C \geq 0.$
\edes
 Assume that these two steps are completed. We now explain how the theorem can be deduced.
From step 2 we deduce that for all $n \geq 0$
\beq
\label{eq:Qnlyap}
Q_{nT_0} V \leq \rho^n V + \frac{C}{1-\rho}.
  \eeq
  Choose $R > \frac{2C}{1-\rho}$ and set $\tilde{K} = \{x \in D: \: V(x) \leq R\}.$ Then $\tilde{K}$ is a compact subset of $D$ and, by step 1, there is  some $n \geq 1$ such that
  \beq
  \label{eq:Qnsmall}
  Q_{nT_0}(x, \cdot) \geq c \nu(\cdot)
  \eeq on $\tilde{K}.$
Now, relying on a version of Harris's theorem proved by Hairer and Mattingly in \cite{Hairer-Mattingly}, (\ref{eq:Qnlyap}) and (\ref{eq:Qnsmall}) imply that  for all $f : D \mapsto \RR$ measurable, and $k \geq 0,$
$$|Q_{knT_0}(f)(x) - \pi(f)| \leq C \gamma^k (1 + V(x)) \|f\|_V$$ for all $x \in D,$ where $0 \leq \gamma < 1$ and $\|f\|_V = \sup_{x \in D} \frac{|f(x)|}{1 + V(x)}.$
Recalling that $h$ is bounded (it lies in $C_0(D)$) and $V = 1/h$, one has $\inf_{x \in D} V(x) > 0.$ This entails the existence of a constant $C > 0$ such that, for all measurable function $f$ such that $\|f\|_V \leq 1$, for all $k \geq 0$, for all $x \in D$,
$$|Q_{knT_0}(f)(x) - \pi(f)| \leq C V(x) \gamma^k.$$
Integrating this last inequality with respect to $\rho(dx)$, for any probability measure $\rho \in \cM_1(D)$,
$$|\rho Q_{knT_0}(f) - \pi(f)| \leq C \rho(V) \gamma^k.$$
Hence, by the property of semigroup for $(Q_t)$, for any $k \in \Z_+$ and $s \in [0,n{T_0}]$, for all $\rho \in \cM_1(D)$,
\begin{align*}
|\rho Q_{s + kn{T_0}}(f) - \pi(f)| 
& \leq C \rho Q_s(V) \gamma^k.\end{align*}
Now, for all $s \in [0,n{T_0}]$ and $\rho \in \cM_1(D)$,
$$\rho Q_s(V) = \rho Q_s(1/h) = \int_D \rho(dx) \frac{e^{\lambda s} P_s^D(\1_D)(x)}{h(x)} \leq e^{\lambda n {T_0}} \rho(1/h).$$
Thus, setting $\alpha=-\ln \gamma$, we deduce that, for any $t \geq 0$ and $f$ such that $\|f\|_{1/h} \leq 1$,
$$|\rho Q_tf - \pi(f)| \leq C e^{\lambda n {T_0}}\rho(1/h)\gamma^{-1} e^{- \alpha t}.$$
In particular, 
by definition of $\pi$, for any $f$ such that $\|f\|_\infty \leq 1$ and $\rho \in \cM_1(D)$ such that $\rho(1/h) < + \infty$, up to a change in the constant $C>0$,
$$\left|\rho Q_t \left[\frac{f}{h}\right] - \frac{\mu(f)}{\mu(h)}\right| \leq C \rho(1/h) e^{-\alpha t}.$$
For any $\rho \in \cM_1(D)$, denote $\rho_h(dx) := \frac{h(x)\rho(dx)}{\rho(h)}\in\mathcal{M}_1(D)$. Remark then that, for all $\rho \in \cM_1(D)$, $\rho_h(\frac{1}{h}) = \frac{1}{\rho(h)}<+\infty$. Then the last inequality (applied to $\rho_h$ instead of $\rho$) entails that, for any $f$ such that $\|f\|_\infty \leq 1$ and $\rho \in \cM_1(D)$, for all $t \geq 0$,
\begin{equation}\label{ref}\frac{\mu(f)}{\mu(h)} - \frac{C}{\rho(h)} e^{-\alpha t} \leq (\rho_h) Q_t \left[\frac{f}{h}\right] \leq \frac{\mu(f)}{\mu(h)} + \frac{C}{\rho(h)} e^{-\alpha t}.\end{equation}
Moreover, using that $P_t^D f(x) = e^{-\lambda t} h(x) Q_t[f/h](x)$,
$$\frac{\rho P_t^D f}{\rho P_t^D 1} = \frac{( \rho_h) Q_t[f/h]}{(\rho_h) Q_t[\1_D/h]}.$$
Thus, fixing $\rho \in \cM_1(D)$ and denoting $t_\rho := \frac{1}{\alpha} \log\left(\frac{C \mu(h)}{\rho(h)}\right)$, this equality and \eqref{ref} entail that, for any $t > t_\rho$,
\begin{equation}
\label{encadrement}\frac{\frac{\mu(f)}{\mu(h)} - \frac{C}{\rho(h)} e^{-\alpha t}}{\frac{1}{\mu(h)} + \frac{C}{\rho(h)} e^{-\alpha t}} \leq \frac{\rho P_t^D f}{\rho P_t^D 1} \leq \frac{\frac{\mu(f)}{\mu(h)} + \frac{C}{\rho(h)} e^{-\alpha t}}{\frac{1}{\mu(h)} - \frac{C}{\rho(h)} e^{-\alpha t}}.
\end{equation}
Therefore, since $|\mu(f)|\leq\|f\|_\infty\leq 1$,
for all 
$t \geq 1 + t_\rho$,
\begin{multline*}
\mu(f) - 2 \frac{C\mu(h)}{\rho(h)} e^{-\alpha t} \leq \frac{\rho P_t^D f}{\rho P_t^D 1} \leq
\mu(f) + 2\frac{\frac{C\mu(h)}{\rho(h)} e^{-\alpha t}}{1 - \frac{C\mu(h)}{\rho(h)} e^{-\alpha t}}\\ \leq
\mu(f) + 2\frac{\frac{C\mu(h)}{\rho(h)} e^{-\alpha t}}{1 - e^{-\alpha}}.
\end{multline*}
Hence, we have proved that, for all
$t > 1 + t_\rho$,
$$\left\| \frac{\rho P_t^D}{\rho P_t^D 1} - \mu \right\|_{TV} \leq \frac{2C}{1-e^{-\alpha}} \frac{\mu(h)}{\rho(h)} e^{-\alpha t}.$$
If
$t \leq 1 + t_\rho$,
$$\left\| \frac{\rho P_t^D}{\rho P_t^D 1} - \mu \right\|_{TV} \leq 2 \leq 2 e^{\alpha (1 + t_\rho)} e^{-\alpha t} = \frac{2 e^\alpha C \mu(h)}{\rho(h)} e^{-\alpha t}.$$
To sum up, there exists two constants $C, \alpha > 0$ such that, for all $\rho \in \cM_1(D)$ and $t \geq 0$,
$$\left\| \frac{\rho P_t^D}{\rho P_t^D 1} - \mu \right\|_{TV} \leq \frac{C \mu(h)}{\rho(h)} e^{-\alpha t}.$$ 
Hence, it remains to prove steps 1 and 2 to conclude the proof.

 We start with a preliminary step, ensuring that there exist a
 non trivial measure $\zeta$ with $\zeta(U) > 0,$ and  positive numbers $T_0 > \eps > 0$ such that for all $x \in U, T_0-\eps \leq t \leq T_0,$
 \begin{align}
 \label{eq:PtDzeta}
 P_t^D(x, \cdot) \geq \zeta(\cdot).
 \end{align}

\underline{Preliminary step:} Since $U$ is $D$-accessible from all $x\in D$, we deduce that
$$
\xi G^D\Ind_U(x)=\int_0^\infty \xi P^D_t\1_U(x)\,dt>0.
$$
Hence there exists $t_\xi>0$ such that $\xi P^D_{t_\xi}\1_U>0$. Setting
$
\zeta'=\xi P^D_{t_\xi}
$
and $T'=T+t_\xi$,
we then obtain $\zeta'(U)>0$ and, for all $x\in U$,
$$
P^D_{T'}(x,\cdot)=\delta_x P^D_TP^D_{t_\xi}\geq \xi P^D_{t_\xi}=\zeta'.
$$
Now, by Fatou Lemma and right continuity of paths (like in the proof of Lemma \ref{lem:GDaccess}),  $\liminf_{t\to 0} \zeta' P^D_t(U)\geq \zeta'(U)>0.$ Then, there exist $\varepsilon,\delta>0$ such that $\zeta' P^D_t\1_U\geq \delta$ for all $t\in[0,\varepsilon]$. Hence, for all $x\in U$ and all $t\in[0,\varepsilon]$,
$$
P^D_{2T'+t}(x,\cdot)\geq \zeta' P^D_{T'+t}\geq \zeta' P^D_t(\1_UP^D_{T'}(\cdot))\geq \delta\zeta'.
$$
Setting $T_0=2T'+\varepsilon$ and $\zeta=\delta\zeta'$, we conclude that~\eqref{eq:PtDzeta} holds true.

\underline{Proof of Step 1:}
 For $\delta > 0$ let  $O_{\delta} = \{x \in D: \: G^D(x,U) > \delta \}.$ Using a similar argument as for
   Proposition~\ref{lem:out}, it is not hard to see  that $O_{\delta}$ is open in $D.$ 
  By (iii), the family $(O_{\delta})_{\delta > 0}$ covers $D.$ Thus, for every compact set $\tilde{K} \subset D$ there exists $\delta > 0$ such that for all
$x \in \tilde{K}, G^D(x,U) \geq \delta.$
Now, relying on Proposition \ref{lem:out}, 
one can choose $S > 0$ large enough so that, for all $x \in \tilde{K},$
$$\int_0^S P_t^D(x,U) dt > \frac{\delta}{2}.$$ Consequently, for all $x \in \tilde{K}$ there is some $0 \leq t_x  \leq S$ such that
$$P_{t_x}^D(x,U) \geq \frac{\delta}{2S} =: \delta'.$$
On the other hand, $h$ is bounded below by a positive constant on $U,$ because for all $x \in U$ $h(x) = e^{\lambda T_0} P_{T_0}^Dh(x) \geq e^{\lambda T_0} \zeta(h) > 0.$
Hence,
$$Q_{t_x}(x,U) \geq \delta''$$ for some $\delta''> 0$ and all $x \in \tilde{K}.$
By~\eqref{eq:PtDzeta}, there exists $c' > 0$ such that
$$Q_t(x, \cdot) \geq c' \zeta'(\cdot)$$ for all $x \in U$ and $T_0-\eps \leq t \leq T_0$, where $\zeta'(f)=\zeta(fh)$.
Choose now $n$ sufficiently large so that $\frac{S}{n} < \eps.$
Then for all $x \in \tilde{K},  n T_0 = t_x + n \tau_x$ for some $\tau_x \in [ T_0-\eps, T_0].$
Thus
$$Q_{nT_0}(x,\cdot) \geq \int_U Q_{t_x}(x,dy) Q_{n \tau_x}(y, \cdot) \geq \delta''  (c' \zeta'(U))^{n-1} c' \zeta'(\cdot).$$

\underline{Proof of Step 2:}
By Markov inequality,
$$\mathbb{P}_x(\tau_D^{out} > T_0) \leq \frac{\mathbb{E}_x(\tau^{out}_D)}{T_0} = \frac{G^D(x,D)}{T_0}.$$ Let $\theta > 0$ be such
that $\rho = e^{\lambda T_0} \theta < 1.$ By the assumption that $G^D(C_b(D)) \subset C_0(D),$  there exists $\eta > 0$ and $C' \geq
0$ such that $$\mathbb{P}_x(\tau_D^{out} > T_0) \leq \theta +  C' \Ind_{\{x \in D :\: d(x, \partial_K D) \geq \eta\}}.$$
Then
$$Q_{T_0}(V)(x) = e^{\lambda T_0} \frac{\mathbb{P}_x(\tau_D^{out} > T_0)}{h(x)} \leq \rho V(x) + \frac{C'e^{\lambda T_0}}{\inf_{y\in
    D,\, d(y,\partial_K D)\geq\eta}h(y)}.$$
\qed

\section{Application to SDEs}
\label{sec:sde}
The main purpose of this section is to prove Theorems \ref{th:mainresultforsde} and~\ref{th:sdehypo}, to complete the proof
  of Corollary~\ref{cor:sde}, and to provide criteria allowing to check their assumptions. {\color{black} In order to do so, we work
  in the SDE settings of Section~\ref{sec:sde0} and use the results of Section~\ref{sec:main} when $M = \RR^n$, $D=\cD$ is a
  connected, bounded, open subset of $\RR^n$ and $K=\overline \cD$.  For simplicity, we denote $\partial {\cal D}$ instead of
    $\partial_K {\cal D}$ in this section.

  Recall that the sets $B(\cD)$, $C_b(\cD)$ and $C_0(\cD)$ denote respectively the space of bounded measurable functions on $\cD$, of
  bounded continuous functions on $\cD$ and  the closure of the set of compactly supported functions on ${\cal D}$
    (i.e.\ the set of bounded continuous functions $f$ on $\cD$ such that $f(x) \rar 0$ when  $d(x,\partial \cD)\to 0${\color{blue})}. These function spaces are equipped with the norm of uniform convergence.

  We let $(P_t)$ denote the Markov semigroup induced by the SDE (\ref{eq:sde}), and by $(P_t^\cD)$ and $G^\cD$ respectively the  sub-Markov semigroup and Green kernel of the sub-Markov semi-group induced by the SDE (\ref{eq:sde}) killed when it exits $\cD$: for all $f\in B(\cD)$, $x\in \cD$ and $t\geq 0$,
  \begin{align*}
  P^\cD_t f(x)=\mathbb E_x\left(f(X_t)\Ind_{t<\tau^\cD_{out}}\right)\text{ and }G^\cD f(x)=\int_0^\infty P_u^\cD f(x)\,du.
  \end{align*}
  }

\subsection{Accessibility}
\label{sec:accessibility}
 In this section, we recall that the general notion of accessibility defined in
  Section~\ref{sec:main} reduces to the one of Section~\ref{sec:sde0} for diffusions and we  provide simple criteria ensuring that Hypothesis (H1) (accessibility of $\mathbb R^n \setminus \bD$ from $\bD$) and condition (ii) of Theorem \ref{th:mainresultforsde} (accessibility of $\cD_{\eps}$ from $\cD$)   are satisfied.

Recall that $y(u,x, \cdot)$ denotes the  maximal solution to the control system (\ref{controlsde}) starting from $x$ (i.e~ $y(u,x,0) = x$).

The following proposition easily follows from the
 celebrated Stroock and Varadhan's support theorem \cite{StrVar72} (see also Theorem 8.1, Chapter VI in \cite{Ikeda}). It justifies
 the terminology used in Definition~\ref{def:SDEaccessible} and makes the link with Definitions~\ref{def:accessible} and~\ref{def:Daccessible}.
 \bprop
 \label{prop:control} Assume that the $S^j, j \geq 1,$ are bounded with bounded first and second derivatives and $S^0$ is Lipschitz and bounded. Let $x \in \mathbb R^n$ and $U \subset \mathbb R^n$
 open. Then $U$ is accessible (respectively $\cD$-accessible) by $(P_t)$ from $x$ (see Definitions \ref{def:accessible} and \ref{def:Daccessible}) if it is accessible (respectively $\cD$-accessible)  by $\{S^0, (S^j)\}$ from $x$ (see Definition \ref{def:SDEaccessible}). \eprop
 \brem {\rm The boundedness assumption is free of charge here, since - by compactness of $\bD$ - we can always
   modify the $S^j$ outside $\bD$ so that they have compact support.} \erem
As an
   illustration of this latter proposition, we provide an  elementary  proof of the following result, originally due to Pinsky \cite{Pinsky}.
   \bprop
\label{cor:pinsky}
Suppose there exists $\tilde{x} \in \mathbb R^n \setminus \bD$ and $\delta > 0$ such that for all $x \in \cD$
\beq
\label{eq:pinsky}\sum_{j = 1}^m \langle S^j(x), x -\tilde{x} \rangle^2 \geq \delta \|x-\tilde{x}\|^2.
\eeq Then $\mathbb R^n \setminus \bD$ is accessible by $\{S^0, (S^j)\}$ (hence by $(P_t)$) from all $x \in \bD.$
\eprop
\prf
Note that by continuity and compactness (of $\bD \times \{x \in \RR^n:
 \: \|x\| = 1\}$), one can always assume that (\ref{eq:pinsky}) holds true on some larger open bounded domain $\cD'$ with $\bD \subset \cD'.$
 For  $x \in \cD$ and $j \in \{1, \ldots, m\}$ set $$u^j(x) =  -  \frac{1}{2 \delta \eps} \langle S^j(x), x -\tilde{x} \rangle$$ where $\eps$ will be chosen later. Consider the ODE
$$\dot{x} = S^0(x) + \sum_j u^j(x) S^j(x)$$ and set $v(t) = \|x(t)-\tilde{x}\|^2.$ Then, as long as $x(t) \in \cD',$
$$\frac{dv(t)}{dt} \leq -\frac{1}{\eps} v(t) + a$$ where $a = \sup_{x \in \cD'} 2 (\langle S^0(x), x-\tilde{x} \rangle).$
 Thus $$v(t) = \|x(t)-\tilde{x}\|^2 \leq e^{-t/\eps} (v_0 - a \eps) + a \eps.$$
   One can choose $\eps$ small enough so that $x(t)$
 meets $\cD' \setminus \bD.$ \qed
The next result (Proposition \ref{cor:inward}) provides a natural condition ensuring that condition $(ii)$ of Theorem \ref{th:mainresultforsde} holds.

Suppose that $\cD$ satisfies the exterior sphere condition  as defined in Definition \ref{def:unitoutward}. Let ${\cal N}_p$ denote the set of unit outward normal vectors at $p \in \partial \cD$  and let
$${\cal N} = \{(p,v):  \: p \in \partial \cD, v \in {\cal N}_p\}.$$
\begin{definition}
  \label{def:point-inward}
We say that a vector $w \in \RR^n$  {\em points inward} $\cD$ at $p \in \partial \cD$ if
$$\langle w, v \rangle \leq 0 \mbox{ for all } v \in {\cal N}_p, \mbox{ \bf and } \langle w, v \rangle < 0 \mbox{ for at least one } v \in {\cal N}_p.$$
 We say that it {\em points strictly inward} $\cD$ at $p$ if  $\langle w, v \rangle < 0$ {\bf for all} $v \in {\cal N}_p.$
We say that a vector field $F$ points inward (respectively strictly inward) $\cD$ if  $F(p)$ points inward (strictly inward) $\cD$ at $p$ for all $p \in \partial \cD.$
\end{definition}

\blem
\label{lem:attractor}
Suppose that $\cD$ satisfies the exterior sphere condition. Let $F$ be a Lipschitz vector field pointing inward $\cD$ and let $\Phi = \{\Phi_t\}$ be its flow. Then

\bdes
\iti $\Phi_t(\bD) \subset  \cD$ for all $t > 0;$
\itii  There exists a compact set $A \subset \cD$ invariant under $\{\Phi_t\}$  (i.e~$\Phi_t(A) = A$ for all $t \in \RR$) such that for all $x \in \cD$ $\omega_{\Phi}(x) \subset A,$  where $\omega_{\Phi}(x)$ stands for the omega limit set of $x$ for $\Phi.$
\edes
\elem
\prf We first show that $\Phi_t(\bD) \subset \bD$ for all $t \geq 0.$ Suppose not. Then, for some $p \in \partial \cD$ and $\eps > 0,$
$d(\Phi_t(p),\bD) > 0$ on $]0,\eps].$ The function $$t \rar V(t):= d(\Phi_t(p),\bD),$$ being Lipschitz on $[0,\eps]$ it is absolutely
continuous, hence almost everywhere derivable and  $V(t) =  \int_0^t \dot{V}(u) du.$  Let $0 < t_0 \leq \eps$ be a point at which it
is derivable, $x_0 = \Phi_{t_0}(p)$ and $p_0 \in  \partial \cD$ be such that $\|x_0-p_0\| =d(x_0,\bD).$ In particular, $\frac{x_0-p_0}{\|x_0-p_0\|}\in\mathcal{N}_{p_0}$.Then for all $t > 0,$
 $$\frac{V(t_0 +t) - V(t_0)}{t} = \frac{d(\Phi_t(x_0),\bD) - d(x_0,\bD)}{t} \leq \frac{\|\Phi_t(x_0)-p_0\| - \|x_0-p_0\|}{t}.$$ Letting
 $t \rar 0,$ using that $F(p_0)$ points inward $\cD$ at $p_0$, we get
 \begin{align*}
  \dot{V}(t_0) \leq \frac{\langle F(x_0), x_0 - p_0\rangle}{\|x_0 - p_0\|} & \leq  \frac{\langle F(x_0) - F(p_0), x_0 -
    p_0\rangle}{\|x_0 - p_0\|} \\ & \leq L \|x_0 - p_0\| = L V(t_0),
 \end{align*}
  where $L$ is a Lipschitz constant for $F$.
  By Gronwall's lemma we then get that for all $0 < s < t \leq \eps,$ $V(t) \leq e^{L (t-s)} V(s).$ Since $V(0) = 0,$ $V$ cannot be
  positive. This proves the desired result. Note that this first result doesn't require that the exterior sphere condition holds
  everywhere
but only that, for all $p \in \partial \cD$ at which it holds,  $\langle F(p),v \rangle \leq 0$ for all $v \in  {\cal N}_p$ (compare with Theorem 2.3 in \cite{Bony}).

  We now show that $\Phi_t(\bD) \subset \cD$ for $t > 0.$ This amounts to show  that for all $p \in \partial \cD$ and   all $\eps>0$ small
  enough, $\Phi_{-\eps}(p) \in \mathbb R^n \setminus \bD.$
So let $p \in \partial \cD$ and choose $v \in {\cal N}_p$ such that $\delta := - \langle F(p),v\rangle > 0.$ By assumption $d(p + rv,\bD) = r$ for some $r > 0.$
Thus for all  $\eps > 0$ such that $\eps \|F(p)\|^2 <   r \delta,$  $$\|(p - \eps F(p)) - (p+ rv)\|^2 < r^2 - \eps \delta r.$$
Since $\|\Phi_{-\eps}(p) - (p -\eps F(p)\| = o(\eps),$ this shows that $\|\Phi_{-\varepsilon}(p)-(p+rv)\|<r$ for
  $\varepsilon>0$ small enough, so that
$\Phi_{-\eps}(p) \not \in \bD$.

Assertion $(ii)$ is a consequence of $(i).$ It suffices to set $$A = \cap_{t \geq 0} \Phi_t(\bD).$$
\qed
For $(p,v) \in {\cal N},$ set
$R(p,v) =  \sup \{r > 0 \: : d(p + rv,\bD) = r\}   \in (0,\infty],$ and  let $$R_{\partial \cD} =\inf_{(p,v) \in {\cal N} } R(p,v).$$
\begin{definition}
\label{def:uniformspherecond}
We say that $\cD$ satisfies the strong exterior sphere condition if it satisfies the exterior sphere condition and  $R_{\partial \cD} \neq 0.$
\end{definition}
\brem
\label{rem:Nicoexamples}
{\rm If $\cD$ is convex or $\partial \cD$ is $C^2,$  then $\cD$ satisfies the strong exterior sphere condition. However,
    the following example shows that the exterior sphere condition and the strong exterior sphere condition are not equivalent:
For $1 \leq \alpha \leq 2,$ let $\cD^{\alpha} \subset \RR^2$ be defined as
$$\cD^{\alpha}= \{(x,y) \in \RR^2 \: : 0 < x < 1, |y| < x^{\alpha}\}.$$  Then $\cD^{\alpha}$ satisfies the exterior sphere condition  for all $1 \leq \alpha \leq 2$ but not the strong one  for $1 < \alpha < 2.$ Indeed, for such an $\alpha,$ ${\cal N}_{(0,0)} = \{v \: : \|v\|= 1, v_1 < 0\}$ and
$$\lim_{v \rar (0,1), v \in {\cal N}_{(0,0)}} R((0,0), v) = 0.$$}
\erem
\bprop
\label{cor:inward}
Assume that $\cD$ satisfies the strong exterior sphere condition. Recall that  ${\cal N}_p$ denotes the set of unit outward normal vectors at $p \in \partial \cD.$
\bdes
\iti For each $p \in \partial \cD,$ the two following conditions 
are equivalent:
\bdes
\ita ${\cal N}_p \cap - {\cal N}_p = \emptyset$ and for each $v \in {\cal N}_p$  $$\langle S^0(p),v\rangle < 0 \mbox{ or } \sum_{i = 1}^m \langle S^i(p),v \rangle^2 \neq 0;$$
\itb There exists a vector $w \in \mathsf{Span} \{S^1(p), \ldots, S^m(p)\}$ such that $S^0(p) + w$ points strictly inward $\cD$ at $p.$
\edes
\itii If for all $p \in \partial \cD$ condition $(i)-(a)$ (or $(i)-(b)$) holds, then  $$\cD_{\eps} = \{x \in \cD: \: d(x,\partial \cD) > \eps\}$$ is
 $\cD$-accessible by $\{S^0, (S^j)\}$ from all $x \in \cD,$ for some $\eps > 0.$
\edes
\eprop
\prf We first prove Point $(ii).$ Assume that $(i)-(b)$ holds at every $p \in \partial \cD$ and prove that $\cD_{\eps}$ is accessible for some $\eps > 0.$ The assumption $R_{\partial \cD} > 0$ makes ${\cal N}$ closed (hence compact). Indeed, if  $(p_n,v_n) \rar (p,v)$ with $(p_n,v_n) \in {\cal N},$ then, for any $0 < r < R_{\partial \cD}$,  $d(p_n + r v_n, \bD) = r.$ Thus $d(p+rv,\bD) = r.$

 This has the consequence that, if a continuous vector field $F$ points strictly inward $\cD$ at $p \in \partial \cD,$ it points strictly inward $\cD$ at $q \in \partial \cD$ for all $q$ in a neighborhood of $p.$ Therefore, by compactness, there exists a covering of $\partial \cD$ by  open sets $U_1, \ldots, U_k,$  and vector fields $W_1, \ldots, W_k \in \mathsf{Span}\{S^1, \ldots, S^m\}$ such that for all $p \in \partial \cD \cap U_i$ $F_i(p) := S^0(p) + W_i(p)$ points strictly inward $\cD$ at $p.$ Set $U_0 = \mathbb R^n \setminus \partial \cD$ and let $\{\rho_i\}_{i = 0,\ldots,k}$ be a partition of unity subordinate to $\{U_i\}_{i = 0,\ldots,k}.$ That is $\rho_i \in C^{\infty}(\mathbb R^n), \rho_i \geq 0, \sum_{i = 0}^k \rho_i = 1,$ and $supp(\rho_i) \subset U_i.$
 Define $F = \rho_0 S^0 + \sum_{i = 1}^k \rho_i F_i.$ Then, $F$ points strictly inward $\cD$  and writes
 \beq
 \label{eq:defF}
 F = S^0 + \sum_{i = 1}^m u^i S^i
  \eeq
  with $u^i \in C^{\infty}(\mathbb R^n).$  In view of (\ref{eq:defF}), 
 Lemma \ref{lem:attractor} proves the result.

  Point $(i).$ We now prove that conditions $(i)-(a)$ and $(i)-(b)$ are equivalent. The implication $(i)-(b) \Rightarrow (i)-(a)$ is straightforward. We focus on the converse implication.
 Let $$Cone({\cal N}_p) = \{t v, t \geq 0, v \in {\cal N}_p\}$$ and $conv({\cal N}_p)$ be the convex hull of ${\cal N}_p.$
  We claim that
  $$conv({\cal N}_p) \subset Cone({\cal N}_p)$$ and $0 \not \in conv({\cal N}_p).$ To prove the first inclusion, it suffices to show
  that $Cone({\cal N}_p)$ is convex. To shorten notation, assume (without loss of generality) that $p = 0.$ Let $x,y \in Cone({\cal
    N}_0)$ and $0 \leq t \leq 1.$ By definition of ${\cal N}_0,$ $Cone(\mathcal{N}_0)=\{z\in\mathbb{R}^n,\exists r>0\text{ s.t.\ }d(rz,\bD)=\|rz\|\}$,
  so there exists $r > 0$ such that $d(rx,\bD) = \|rx\|$ and $d(ry,\bD) = \|ry\|.$ Thus for all $z \in \bD$
  $$\|rx -z\|^2-\|rx\|^2 = \|z\|^2-2 \langle rx,z \rangle \geq 0.$$ Similarly $\|z\|^2-2 \langle ry,z \rangle \geq 0.$
  Thus
  $$\|r(tx +(1-t)y) -z\|^2 - \|r(tx +(1-t)y)\|^2$$
  $$ = t(\|z\|^2-2 \langle rx,z \rangle) + (1-t)(\|z\|^2 -2  \langle ry,z \rangle) \geq 0.$$
  This proves that $tx +(1-t)y \in Cone({\cal N}_0),$ hence convexity of $Cone({\cal N}_0).$

  The fact that $0 \not \in conv({\cal N}_p)$ follows from the fact that ${\cal N}_p \cap - {\cal N}_p = \emptyset.$ 
  Indeed, suppose to the contrary that $0 = \sum_{i = 1}^k t_i x_i$ with $k \geq 2, x_i \in {\cal N}_p, t_i > 0$ and $\sum_{i = 1}^k t_i = 1.$ Then
  $$- \frac{t_1}{1-t_1} x_1 \in conv(x_2, \ldots, x_k) \subset conv({\cal N}_p) \subset cone({\cal N}_p).$$ Thus $-x_1 \in {\cal N}_p.$ A contradiction.

  We shall now deduce the implication $(i)-(a) \Rightarrow (i)-(b)$ from the Minimax theorem (see e.g~ \cite{Sorin}).
  For all $j \in \{1, \ldots, m\}$ set $S^{-j} = - S^j.$ Let $J = \{-m, \ldots, 0, \ldots m\}$ and $$\Delta(J) = \{\alpha \in \RR^J : \: \alpha_j \geq 0, \sum_{j \in J} \alpha_j = 1\}.$$
  By condition $(i)-(a)$ and compactness of ${\cal N}_p$ there exists $\delta > 0$ such that
  for all $v \in {\cal N}_p$ $$\min_{j \in J} \langle S^j(p),v\rangle \leq -\delta.$$ Thus for all $v \in Cone({\cal N}_p)$
  $$\min_{j \in J} \langle S^j(p),v\rangle \leq -\delta \|v\|.$$ The set ${\cal N}_p$ being compact (in finite dimension) its convex hull is also compact by Carath\'eodory's theorem. Thus, because $0 \not \in conv({\cal N}_p),$ $\|v\| \geq \frac{\delta'}{\delta}$ for some $\delta' > 0$ and all $v \in conv({\cal N}_p).$ It then follows that
  $$\sup_{v \in conv({\cal N}_p)} \inf_{\alpha \in \Delta(J)} \langle \sum_{j \in J} \alpha_j S^j(p),v\rangle \leq \sup_{v \in conv({\cal N}_p)} \min_{j \in J} \langle S^j(p),v\rangle \leq -\delta'.$$ 
  By the Minimax theorem, the left hand side also equals
  $$ \inf_{\alpha \in \Delta(J)} \sup_{v \in conv({\cal N}_p)} \langle \sum_{j \in J} \alpha_j S^j(p),v\rangle$$ and this infimum is achieved for some $\beta \in  \Delta(J).$ If $\beta_0 \neq 0$ this implies that
  $$\sup_{v \in conv({\cal N}_p)} \langle S^0(p) + \sum_{j \in J, j \neq 0} \frac{\beta_j}{\beta_0} S^j(p),v\rangle  \leq -\frac{\delta'}{\beta_0} < 0.$$
  If $\beta_0 = 0,$ for $R > 0$ sufficiently large $$\sup_{v \in conv({\cal N}_p)} \langle S^0(p) + R \sum_{j \in J} \beta_j S^j(p),v\rangle  \leq -R \delta' + \|S^0(p)\| < 0.$$
  This concludes the proof.
  \qed
  \brem
   {\rm It follows from  Proposition \ref{cor:inward}  that whenever $\partial \cD$ is $C^2,$ Hypothesis (H2') implies  condition
     $(ii)$ of Theorem \ref{th:mainresultforsde} because at each point $p \in \partial \cD$ there is a unique outward unit normal.
     The  following example shows that this is not true in  general when (H2') is satisfied but $\partial\cD$ is not $C^2$. Let $\cD^1$ be as in Remark \ref{rem:Nicoexamples}, with $\alpha = 1,$
 and let $(X_t)$ be solution to $$dX_t =S^1(X_t) \circ dB_t$$ where $S^1(x,y) = (1,2).$
   At each point $p \in \partial \cD$ there is at least one $v  \in {\cal N}_p$ such that $\langle S^1(p), v \rangle \neq 0$ so that
   condition $(i)-(a)$ is satisfied.  However,$\cD$ does not satisfy the strong exterior sphere condition and, for $0 < \eta < \eps, \cD^1_{\eps}$ is not $\cD$-accessible from  $(\eta,0).$

Observe also that  none of the conditions required in  Proposition \ref{cor:inward} is necessary for $\cD_{\eps}$ to be accessible. Let $\cD^{\alpha}$ be as in Remark \ref{rem:Nicoexamples}, with $1 < \alpha \leq 2$ and let $$dX_t = e_1 \circ dB_t^1 + e_2 \circ dB_t^2,$$ with $(e_1,e_2)$ the canonical basis of $\RR^2.$ As shown in Corollary \ref{cor:sde} (and its proof), $\cD_{\eps}^{\alpha}$ is accessible; while for $1 < \alpha < 2, R_{\partial \cD} = 0$ and for $\alpha = 2, {\cal N}_{0,0} \cap - {\cal N}_{0,0} \neq \emptyset.$}
   \erem
 \subsection{Proof of Theorem \ref{th:mainresultforsde}}
\label{sec:pf-th1.5}

{\color{black}In this Section, we start with general path properties in Lemma~\ref{lem:pathproperties} and show that, under Assumptions~(H1) and~(H2),  $G^\cD(C_b(\cD))\subset C_0(\cD)$ in Lemma~\ref{lem:fellerproperties}. We then prove  Theorem~\ref{th:mainresultforsde}. We conclude this section with Proposition~\ref{prop:H2'impliesH2}, which shows that  Assumption~(H2') implies Assumption~(H2).}

Let $C(\Rp, \mathbb R^n)$ be the set of continuous paths $\eta : \Rp \mapsto \mathbb R^n$ equipped with the topology of uniform convergence on compact intervals.

 Let $(\eta_n)_{n \geq 0}$ be a sequence converging to $\eta$ in $C(\Rp, \mathbb R^n).$
 Set $$\tau^{n,out}_\cD = \inf \{t \geq 0: \: \eta_n(t) \in \mathbb R^n \setminus \cD\}, \tau^{out}_\cD = \inf \{t \geq 0: \: \eta(t) \in \mathbb R^n \setminus \cD\}$$
 Define $\tau^{n,out}_\bD$ and $\tau^{out}_\bD$ similarly. {\color{black} In the following lemmas, $C_0(\bD)$ is the space of bounded continuous functions on $\bD$ vanishing on $\partial \cD$.}
 \blem
 \label{lem:pathproperties}
 \bdes
 \iti
 Suppose $\eta(0) \in \cD.$ Then for all $f  \in C_b(\cD)$ such that $f \geq 0$ and all $t \geq 0,$
 $$\liminf_{n \rar \infty} f(\eta_n(t)) \Ind_{\{\tau^{n, out}_\cD > t\}} \geq f(\eta(t))  \Ind_{\{\tau^{out}_\cD > t\}}.$$ In particular
 $$\liminf_{n \rar \infty} \tau^{n,out}_\cD \geq \tau^{out}_\cD.$$
 \itii
 Suppose $\eta(0) \in \bD.$ Then $$\limsup_{n \rar \infty}  \tau^{n,out}_\bD \leq \tau^{out}_\bD$$ and for all $f \in C_0(\bD), f \geq 0$ and all $t \geq 0,$
 $$\limsup_{n \rar \infty} f(\eta_n(t)) \Ind_{\{\tau^{n,out}_\bD > t\}} \leq f(\eta(t))  \Ind_{\{\tau^{out}_\bD > t\}}.$$
 \edes
 \elem
 \prf
 $(i)$  If $\tau^{out}_\cD \leq t$ the statement is obvious. If $\tau^{out}_\cD > t,$ then $\eta([0,t]) \subset \cD$ so that, for $n$ large enough, $\eta_n([0,t]) \subset \cD.$ That is $\tau^{n,out}_\cD > t$ and the statement follows. The assertion that $\liminf_{n \rar \infty} \tau^{n,out}_\cD \geq \tau^{out}_\cD$ follows by choosing $f = \Ind_\cD.$

 $(ii)$ Suppose to the contrary that $\tau^{n,out}_\bD > \tau^{out}_\bD + \eps$ for some $\eps > 0$ and infinitely many $n.$ Then
 $\eta^n([0,\tau^{out}_\bD + \eps])  \subset \bD.$ Hence $\eta([0,\tau^{out}_\bD + \eps]) \subset \bD.$ A contradiction. The last
 assertion directly follows for $t \neq \tau^{out}_\bD$ and $f \in C_b(\bD)$ {\color{black}(where $C_b(\bD)$ is the set of bounded continuous functions on $\bD$, hence it contains $C_0(\bD)$).} If now $t = \tau^{out}_\bD$ and $f \in C_0(\bD),$ $f(\eta^n(t)) \rar f(\eta(t)) = 0.$
 \qed
The next lemma shows that under (H1) and (H2),  $G^\cD(C_b(\cD))\subset C_0(\cD)$. In this lemma, $P_t^\bD$ denotes the semigroup defined as $P_t^\bD f(x) = \mathbb{E}_x(f(X_t) \Ind_{\tau_\bD^{out} > t})$ for all $f$ bounded and measurable on $\bD$.
\blem
\label{lem:fellerproperties}
Suppose that Hypothesis (H1) holds. Then,
 \bdes
 \iti
 For all $f \geq 0, f \in C_b(\cD), x \in \cD$ and $t \geq 0,$
 $$\liminf_{y \rar x} P_t^\cD f(y) \geq P_t^\cD f(x),$$
 \itii
 For all $f \geq 0, f \in C_0(\bD), x \in \bD$ and $t \geq 0$
 $$\limsup_{y \rar x, y \in \bD} P_t^\bD f(y) \leq P^\bD_t f (x)$$
 \itiii Suppose, in addition, that Hypothesis (H2) holds.  Then,
 \bdes
 \ita For all  $x \in \cD,$ $$\mathbb{P}_x(\tau_\cD^{out} = \tau_\bD^{out}) = 1;$$
 \itb For all $f \in C_0(\cD)$ and $t \geq 0,$
 $P_t^\cD f \in C_b(\cD);$
 \itc For all $f \in C_b(\cD)$ and $t \geq 0,$ $G^\cD f  \in C_0(\cD).$ 
 \edes
 \edes
 \elem
 \prf
 Let $(X_t^x)$ be the (strong) solution to (\ref{eq:sde}) with initial condition $X_0^x = x$. We can always assume that $(X_t^x)_{t \geq 0, \, x \in \mathbb R^n}$ is defined on  the Wiener space  space $C(\Rp, \RR^m)$ equipped with its Borel sigma field and the Wiener measure $\mathbb{P}$ (the law of $(B^1, \ldots, B^m)$). That is $\mathbb{P}_x( \cdot) = \mathbb{P}(X^x \in \cdot).$ Also, for all $\omega \in $ $C(\Rp, \RR^m)$, the map $x \in \mathbb R^n  \rar X^x(\omega) \in C(\Rp, \mathbb R^n)$ is continuous (see for instance \cite{Legall2}, Theorem 8.5).

Assertions $(i)$ and $(ii)$ then follow from Lemma \ref{lem:pathproperties} and Fatou's lemma.

We now pass to the proof of $(iii).$

$(a)$ follows from the strong Markov property,
valid since $(P_t)_{t\geq 0}$ is
  Feller~(see e.g.~\cite{Legall2}, Theorem 6.17), as follows. For all $x \in \cD,$
$$\mathbb{P}_x(\tau_\cD^{out} = \tau_\bD^{out}) = \mathbb{E}_x( \mathbb{P}_{X_{\tau^{out}_\cD}}(\tau_\bD^{out} = 0)) = 1.$$

$(b)$ Let $x \in \cD.$ The property  $\mathbb{P}_x(\tau^{out}_\cD = \tau^{out}_\bD) = 1$ implies that $P_t^\cD f(x) = P_t^\bD
\hat{f}(x)$ for all  $f \in C_b(\cD),$ where $\hat{f}(x)=f(x)$ if $x\in \cD$ and $\hat{f}(x)=0$ if
  $x\in \bD\setminus \cD.$ Hence, by $(i)$ and $(ii),$ $$\lim_{y \rar x} P_t^\cD f(y) = P_t^\cD f(x)$$ for all $f \in C_0(\cD), f \geq 0.$ If now $f \in C_0(\cD)$ it suffices to write $f = f^+-f^{-}$ with $f^+ = \max(f,0)$ and  $f^{-} = (-f)^+.$

 $(c)$ Write $\tau^{x,out}_\cD$  for $\inf \{t \geq 0: \: X_t^x \not \in \cD \}.$ Again, the  property  $\mathbb{P}_x(\tau^{out}_\cD = \tau^{out}_\bD) = 1$ combined with  Lemma \ref{lem:pathproperties}, imply that, almost surely, the maps $x \in \cD \rar \tau^{x,out}_\cD,$ and
$x \in \cD \rar \int_0^{\tau^{x,out}_\cD} f(X_s^x) ds$ are continuous for all $f \in C_b(\cD).$
Also, $$\sup_{x \in \cD} \mathbb{E}\left[ \left(\int_0^{\tau^{x,out}_\cD} f(X_s^x) ds\right)^2\right] \leq \|f\|^2 \sup_{x \in \cD} \mathbb{E}_x[(\tau^{out}_\cD)^2] < \infty$$
where the last inequality follows from Proposition \ref{lem:out}. This shows that the family $(\int_0^{\tau^{x,out}_\cD} f(X_s^x) ds)_{x \in \cD}$ is uniformly integrable. Thus, $x \in \cD \rar \mathbb{E}(\int_0^{\tau^{x,out}_\cD} f(X_s^x) ds) = G^\cD f(x)$ is continuous.

To conclude, observe that $|G^\cD f(x)| \leq \|f\| G^\cD \Ind_\cD(x)$ and that $G^\cD \Ind_\cD(x) = \mathbb{E}(\tau_\cD^{x,out}).$ {\color{black} According to (iii)-(a) of Lemma~\ref{lem:fellerproperties}, we have  $\tau^{x,out}_{\cD}=\tau^{x,out}_{\bD}$ almost surely, and hence, by (ii) of Lemma~\ref{lem:pathproperties},  $\tau^{x,out}_{\cD}\to \tau^{p,out}_{\bD}$. Since $\tau^{p,out}_{\bD}=0$ almost surely by Assumption~(H2) and since $(\tau^{x,out}_{\cD})_{x\in\cD}$ is uniformly integrable by Proposition~\ref{lem:out}, we deduce that  $\mathbb{E}(\tau_\cD^{x,out})$
 converges to  $0$ when $x\to p$.}
 \qed

{\color{black}
    \noindent\textbf{Proof of Theorem~\ref{th:mainresultforsde}: } We assume that
  \bdes
  \iti Hypotheses (H1) and (H2) hold true;
  \itii For some $\eps > 0,$ the set $\cD_{\eps} = \{x \in \cD : \: d(x,\partial \cD) > \eps\}$ is $\cD$-accessible  by $\{S^0, (S^j)\}$ from all $x \in \cD \setminus \cD_{\eps}.$
\edes
    In order to prove Theorem \ref{th:mainresultforsde}, we show that the assumptions  of Theorem \ref{th:main} are satisfied  with $M=\mathbb R^n$, $D=\cD$, $K=\bD$, $X$ following the dynamic defined by~\eqref{eq:sde}, $P^D=P^\cD$ and $G^D=G^\cD$. According to Lemma~\ref{lem:fellerproperties} $(iii)-(c)$, Hypothesis~(H1) entails that $G^\cD(C_b(\cD)) \subset C_0(\cD).$ In addition, Condition $(ii)$ of Theorem~\ref{th:mainresultforsde} and Proposition~\ref{prop:control} entail that Assumption $(ii)$ of Theorem~\ref{th:main} holds true. In particular, we deduce that Theorem~\ref{th:main} applies and hence that there exists a QSD for $(X_t)$ on $\cD$.
 \qed }

In the next result, we prove that (H2') implies (H2).
\bprop
\label{prop:H2'impliesH2}
Let $p \in \partial \cD.$ Assume that there exist a unit outward vector $v$ at $p$ 
 and
$j \in \{1, \ldots, m\}$ such that $\langle  v, S^j(p) \rangle \neq 0.$ Then, $p$ is regular with respect to $\mathbb R^n \setminus \bD$ for (\ref{eq:sde}). In particular, Hypothesis $(H2')$ implies $(H2).$
\eprop
\prf
Without loss of generality we can assume that $j = 1.$
 By Definition \ref{def:unitoutward}, there is $r > 0$ such that $d(p+rv,\bD) = r.$ Let $\Psi:  \: \RR^n \mapsto \RR$ be defined as $$\Psi(x) = r^2 - \|x -(p + rv)\|^2.$$ Then $$\Psi(x) > 0 \Rightarrow x \in \mathbb R^n \setminus \bD,$$ and
 $\langle \nabla \Psi(p), S^1(p) \rangle \neq 0.$  Hence,  for some  neighborhood $U$  of $p,$ $$(\langle \nabla \Psi(x), S^1(x) \rangle)^2 \geq a >  0$$ on $\overline{U}.$

 By Ito's formula, since $\Psi(p)= 0,$
 $$\Psi(X^p_{t \wedge \tau^{out}_U})= \int_0^{t \wedge \tau^{out}_U} L \Psi(X_s^p) ds + M_{t \wedge \tau^{out}_U}$$
 where
  $$L \Psi = S^0(\Psi) + \frac{1}{2} \sum_{j = 1}^m (S^j)^2 (\Psi),$$
 $$M_t = \sum_{j = 1}^m \int_0^t  \sigma^j(X_s^p) dB^j_s,$$
 and
 $\sigma^j(x)$ is any bounded measurable function coinciding with
 $\langle \nabla \Psi(x), S^j(x) \rangle$ on $\overline{U}.$ In the definition of $L$ above, we used the standard notation for
 differential operators defined from vector fields: $S^i(f)(x)=\langle S^i(x), \nabla f(x)\rangle$. For convenience, we set $\sigma^1(x) = \sqrt{a}$ for $x \not \in \overline{U}.$
Therefore,
\beq
\label{eq:lowerpsi}
\Psi(X^p_{t \wedge \tau^{out}_U}) \geq   M_{t \wedge \tau^{out}_U} - b (t \wedge \tau^{out}_U)
\eeq
with $b = \sup_{x \in \overline{U}} |L\Psi(x)|,$
and
$$\langle M \rangle_t = \sum_{j = 1}^m \int_0^t \sigma^j(X_s^p)^2 ds \geq at.$$
By Dubins-Schwarz Theorem (see e.g~ \cite{Legall2}, Theorem 5.13) there exists a Brownian motion $(\beta)$ such that
$M_t = \beta_{\langle M \rangle_t}$ for all $t \geq 0.$
Thus, for all $\eps > 0$
$$\sup_{0 \leq t \leq \eps} (M_t - b t) \geq \sup_{0 \leq t \leq \eps} (\beta_{\langle M \rangle_t} -\frac{b}{a} \langle M \rangle_t)
 \geq \sup_{0 \leq t \leq a \eps} (\beta_t - \frac{b}{a}t). $$
 By the law of the iterated logarithm for the Brownian motion, the right hand term is almost surely positive.

Now, using   (\ref{eq:lowerpsi}), it  follows that $\tau^{out}_\bD \leq \eps$
 almost surely on the event $\{\tau^{out}_U > \eps \}.$ Thus $\mathbb{P}_p(\tau^{out}_\bD = 0) = 1$  because $\mathbb{P}_p(\tau^{out}_U > 0) = 1.$
 \qed

 \subsection{Proof of Theorem \ref{th:sdehypo}}
 \subsubsection*{Consequences of Hörmander conditions and Bony's results~\cite{Bony}}
 We assume throughout all this subsection, as in Theorem~\ref{th:sdehypo}, that: \bdes
 \iti Hypotheses (H1) and (H2') holds true;
 \itii The weak
 H\"ormander condition (see Definition \ref{def:hormanderconditions}) is satisfied at every $x \in \bD.$
 \edes

 \blem
 \label{lem:hypocomp}
 The operator $G^\cD$ is a compact operator on $C_0(\cD).$
 \elem
 \prf
 By Lemma \ref{lem:fellerproperties}, (iii) (c) and Proposition \ref{lem:out}, the operator $G^\cD$ is a bounded operator on
 $C_b(\cD)$ whose image is in $C_0(\cD).$ It then defines a bounded operator on $C_0(\cD).$ To prove compactness, we rely on the
 following (easy) consequence of Ascoli's theorem (see e.g~\cite{Munkres}, exercise 6, chapter 7): A family
 ${\cal F} \subset C_0(\cD)$ is relatively compact if (and only if) it is bounded, equicontinuous
 at every point $x \in \cD$ and {\em vanishes uniformly}
 at $\partial \cD$ \bla{(regardless of the smoothness of $\partial \cD$)}; that is for all $\eps > 0$ there exists a compact set $W \subset \cD$ such that $|f(x)| \leq \eps$ for all
 $f \in {\cal F}$ and $x \in \cD \setminus W.$

 Let ${\cal F} = \{G^\cD(f):  \: f \in C_0(\cD), \|f\| \leq 1\}.$ Then ${\cal F}$ is bounded and   vanishes uniformly  at $\partial \cD$ because for all $f$ in the unit ball of $C_0(\cD)$  $|G^\cD(f)| \leq G^\cD \Ind_\cD \in C_0(\cD).$

 We now prove  the equicontinuity property.
 Let $L$ be the differential operator defined as $L= S^0 + \frac{1}{2} \sum_{j = 1}^m (S^j)^2.$

We claim that for every $f \in C_0(\cD),$ the map $g = G^\cD f$ satisfies $L g  = - f$ on $\cD,$ in the sense of distributions. We proceed in two steps: first we will show how the result can be proved from the claim and second we will prove the claim.

\smallskip\noindent \textit{Step 1.}
First, assume that the claim is proved. Then, by  a theorem due to  Rothschild and Stein (Theorem 18 in \cite{rothschild1976}) there exists
$0 < \alpha < 1$ depending only on the family $S^0, S^1, \ldots, S^m$ such that  $G^\cD(C_0(\cD)) \subset \Lambda^{\alpha}(\cD).$
 Here (see \cite{rothschild1976}, \bla{p.\ 301 and \cite{Stein1970}, p.\ 141}),  $\Lambda^{\alpha}(\cD)$ denotes the set of \bla{locally
 $\alpha$-H\"older functions, that is the set of} $f : \cD \mapsto \RR$ such that for every compact $W \subset \cD,$
\[
\|f\|_{\alpha,W}:=\sup_{x\in W} |f(x)|+\sup_{x\neq y\in W}|f(x)-f(y)|/|x-y|^\alpha<\infty.
\]

Let now $x_0 \in \cD,$ $W \subset \cD$ be a compact neighborhood of $x_0,$ $C(W)$  the space of continuous functions on $W$ equipped
with the uniform norm, $C^{\alpha}(W)$ the Banach space of $\alpha$-H\"older functions on $W$ equipped with the corresponding Hölder norm
\bla{$\|.\|_{\alpha,W}$} and
$i_W : C_0(\cD) \mapsto C(W), f \mapsto f|_W.$ The operator $i_W \circ G^\cD : C_0(\cD) \mapsto C(W)$ is bounded,
  where $\circ$ denotes here the composition of functions. Its graph is then closed in $C_0(\cD) \times  C(W)$ and consequently also in $C_0(\cD) \times C^{\alpha}(W).$ Hence, by the closed graph theorem,  it  is a bounded operator from $C_0(\cD)$ into  $C^{\alpha}(W).$ In particular, for all $f$ in the unit ball of $C_0(\cD)$ and $x,y \in W$ $|G^\cD(f)(x) - G^\cD f(y)| \leq \kappa \|x-y\|^{\alpha}$ for some $\kappa$ depending only on $W.$ This proves that ${\cal F}$ is equicontinuous at  $x_0.$

\smallskip\noindent\textit{Step 2.}  We now pass to the proof of the claim.

 By a theorem of Bony (\cite{Bony}, Th\'eor\`eme 5.2), for any $a > 0$  and $f \in C_b(\bD)$ {\color{black} (the set of bounded continuous functions on $\bD$)}, the Dirichlet problem:

  \beq
  \label{eq:diricha}
  \left \{ \begin{array}{l}
  L g - a g  = - f \mbox{ on } \cD \mbox{ (in the sense of distributions)} \\
  g|_{\partial \cD}  =  0;
  \end{array} \right.
   \eeq
   has  a unique solution, call it $g_a,$ continuous on $\bD.$ Furthermore, if $f$ is smooth ($C^\infty$) on $\cD$ so is $g_a.$
   Note that the assumptions required for this theorem are implied by \bla{
     the properties of $\partial \cD$ and of the vector fields $S^j$ near $\partial \cD$ given in Assumption~(H2')} and the weak H\"ormander
   condition \bla{assumed in (ii)}.

  Suppose that $f$ is smooth on $\cD$, meaning that $f\in C_b(\bD)\cap C^\infty(\cD)$ (so that $g_a$ is smooth on $\cD$). Then by Ito's formula  $$\left(e^{-{a t \wedge \tau^{out}_\cD}} g_a(X_{t \wedge \tau^{out}_\cD}) +\int_0^{t \wedge \tau^{out}_\cD} e^{-a s} f(X_s) ds\right)_{t \geq 0}$$ is a local martingale. Being bounded, it is  a uniformly integrable martingale. Thus, taking the expectation and letting $t \rar \infty$, we get that
  \beq
  \label{eq:ga}
  \int_0^{\infty} e^{-as} P_s^\cD f(x) ds = \mathbb{E}_x\left(\int_0^{\tau_\cD^{out}} e^{- a s} f(X_s) ds\right) = g_a(x).
  \eeq
In particular, $G^\cD f(x)=\lim_{a\rightarrow 0}g_a(x)$, where the convergence is uniform by
Proposition~\ref{lem:out}.

For every smooth test function $\Phi$ with compact support in $\cD$
$$\langle g_a, L^* \Phi \rangle - a \langle g_a, \Phi \rangle = -\langle f,\Phi \rangle$$
where $\langle h, \Phi \rangle = \int h(x) \Phi(x) dx.$ Letting $a \rar 0$, we get that $\langle G^\cD f,L^*\Phi\rangle=-\langle
f,\Phi \rangle$, that is
$$L G^\cD(f) = -f,$$ in the sense of distributions.
This proves the claim whenever $f$ is smooth on $\cD.$ 
If now $f\in\C_0(\cD)$, let $(f_n)_{n\in\mathbb N}$ be smooth with compact support in $\mathcal D$ with $f_n \rar f$ uniformly on $\bD$. Then $G^\cD(f_n)
\rar G^\cD(f)$ uniformly; and, by the same argument as above, $L G^\cD f= -f$ in the sense of distributions.
\qed

The proofs of the next two lemmas are similar to the proof of Corollary 5.4  in \cite{Ben18}. For convenience we provide details.
\blem
\label{lem:hypopetiteset}
Let $p \in \cD$ be such that $S^i(p) \neq 0$ for some $i \in \{1, \ldots, m\}.$ Then, there exist disjoint open sets $U,V \subset \cD$
 with $p \in U$ and a nontrivial measure $\xi$ on $V$ (i.e.\ s.t.\  $\xi(V) > 0, \xi(\mathbb R^n \setminus V) = 0$) such that for all $x \in U,$
$$G^\cD(x, \cdot) \geq \xi(\cdot).$$
\elem
\prf
 We first assume that  for all $x \in \cD$ there is some $i \in\{1, \ldots, m\}$ such that $S^i(x) \neq 0.$ Then,  by Theorem 6.1 in \cite{Bony},
 for all $a > 0,$  there exists a map $K_a : \bD^2 \mapsto \Rp$  smooth on $\cD^2 \setminus \{(x,x):
 \: x \in \cD\}$ such that for all $f \in C_b(\bD),$ and $x \in \cD$
$$\int_0^{\infty} e^{-as} P_s^\cD f(x) ds  = \int_\cD K_a(x,y) f(y) dy.$$
To be more precise, Theorem 6.1 in  \cite{Bony} asserts that the solution to the
 Dirichlet problem (\ref{eq:diricha}) can be written under the form given by the right hand side of the last equality.
 On the other hand,  we have shown in the proof of Lemma \ref{lem:hypocomp}, that this solution is given by the left hand side.

Fix $a = 1$ and choose $q \neq p$ such that $K_1(p,q) > 0.$ Such a $q$ exists, for otherwise, we would have $\int_0^{\infty} e^{-s} P_s^\cD \Ind_\cD(p) ds = 0.$ That
 is $\tau^{out}_\cD = 0, \mathbb{P}_p$ almost surely.
By continuity of $K_1$ off the diagonal, there exist disjoint neighborhoods $U,V$  of $p$ and $q$ and some $c > 0$ such that
$K_1(x,y) \geq c$ for all $x \in U, y \in V.$ Thus, for all $x \in U,$
$$G^\cD(x, \cdot) \geq \int_0^{\infty} e^{-s} P_s^\cD(x, \cdot) \geq  c Leb(V \cap \cdot)$$ where $Leb$ stands for the Lebesgue measure on $\RR^n.$

We now pass to the proof of the Lemma. Using a local chart around $p$ we can  assume without loss of generality that
 $p = 0, S^1(0) \neq 0$ and $\frac{S^1(0)}{\|S^1(0)\|} = e_1,$ where $(e_1, \ldots, e_n)$ stands for the canonical basis on $\RR^n.$
  Let $B_{\eps} = \{x \in \RR^n: \: \sum_{i = 1}^n |x_i| < \eps \}.$ For $\eps > 0$ small enough $S^1(x) \neq 0$ for all $x \in B_{\eps} $ and
for all $x \in \overline{B_{\eps}}$ there is a vector $u$ normal to $\overline{B_{\eps}}
\setminus B_{\eps}$  such that
$\langle S^1(x), u \rangle \neq 0.$ The preceding reasoning can then be applied with $B_{\eps}$ in place of $\cD.$
Thus, for some disjoint open sets $U,V \subset \B_{\eps}$ with $p \in U$ and for all $x \in U,$
$$G^\cD(x, \cdot) \geq G^{B_{\eps}}(x,\cdot) \geq c Leb(V \cap \cdot).$$
\qed
\blem
\label{lem:hyposmallset}
 Let $p \in \cD$ be such that the  Hörmander condition holds at $p.$
 Then, there exist a  neighborhood $U$ of $p,$ a non trivial measure $\xi$ on $\cD$ and   $T > 0$ such that for all $x \in U,
 P_T^\cD(x,\cdot) \geq \xi(\cdot).$
\elem
\prf Let $V \subset \cD$ be a  neighborhood of $p$ such   that the Hörmander condition holds at every $x \in V.$ Let $P_t^V(x,\cdot)$ stands for the law of the stopped process $X_{t \wedge \tau_V^{out}}^x.$    By (\cite{FoldesHerzog}, Theorem 4.3 (ii)), there exists a nonnegative map $p_t^V(x,y)$ smooth in the variables $(t, x,y) \in \Rp^* \times V \times V$ such that  $$P_t^V(x,A) = \int_A p_t^V(x,y) dy$$ for every Borel subset $A$ of $V.$  Choose $q \in V$ and $T> 0$ such that $p_T^V(p,q):= c  >  0.$ Then, by continuity,  there exist  neighborhoods $U,W \subset V$ of $p$ and $q$ such that $p_T^V(x,y) \geq c/2 > 0$ for all $(x,y) \in  U \times W.$ Thus, for all $x \in U,$
$$P_T^\cD(x,\cdot) \geq P_T^\cD(x, \cdot \cap W) \geq  P_T^V(x, \cdot \cap W) \geq \frac{c}{2} Leb(W \cap \cdot).$$
\qed
Note that the continuity properties of the kernels $K_a$ and  $p_t^V$ used in the proof of the previous lemmas \ref{lem:hypopetiteset} and \ref{lem:hyposmallset}  could also be derived from the results obtained by P. Cattiaux (\cite{Cattiaux}, Theorems 4.14 and 3.35)  via Malliavin's calculus.

\bla{We conclude this section with the statement and proof of Lemma~\ref{lem:hyposmooth}, used to complete the proofs of
  Theorem~\ref{th:sdehypo} (given below) and Corollary~\ref{cor:sde} (given in Section~\ref{sec:sde0}).}
\blem
\label{lem:hyposmooth}
Let $\mu$ be a QSD for $(X_t)$ on $\cD.$ Then $\mu$ has a smooth density with respect to the Lebesgue measure on $\cD.$
 If furthermore the strong Hörmander conditions holds at every $x \in \bD,$ this density is positive. \elem
\prf
Let $\Phi$ be a smooth function with compact support in $\cD.$ By Ito's formula
$(\Phi(X_{t \wedge \tau^{out}_\cD}) - \int_0^{t \wedge \tau^{out}_\cD} L \Phi(X_s) ds)_{t\geq 0}$ is a bounded local
martingale, 
hence a true martingale. Thus, taking the expectation, letting $t \rar \infty$ and using Proposition~\ref{lem:out},
it comes that $G^\cD (L \Phi)(x) = - \Phi(x)$ for all $x \in \cD.$ Let $\mu$ be a QSD with rate $\lambda.$ Then $-\mu (\Phi) = \mu G^\cD (L
\Phi) = \frac{1}{\lambda} \mu (L \Phi).$ This shows that $L^* \mu + \lambda \mu = 0$ on $\cD$ in the sense of distributions. Now,
for all distribution $f$ on $\cD$,
$$L^* f = \tilde{S}^0 f + \frac{1}{2} \sum_{j = 1}^m (S^j)^2 f + T f$$ where $T$ is a smooth function and $S^0 + \tilde{S}^0 \in
Span(S^1, \ldots, S^m).$ Therefore  $L^*$ satisfies the weak H\"ormander property. By H\"ormander Theorem \cite{Hormander}, it is
hypoelliptic. This  implies that $\mu$ has a smooth density.

For $a > 0,$ set $L^*_a f = L^*f - af$ and choose $a$ sufficiently large so that  $L^*_a 1 = T-a < 0.$
If the strong H\"ormander condition is satisfied at every point $x \in \bD,$ the same is true for $L^*_a.$  Now, $L^*_a (-\mu) = (\lambda + a)
\mu \geq 0.$ Therefore, by application of Bony's maximum principle (\cite{Bony}, Corollary 3.1),   if the density of $-\mu$  vanishes
at some $x \in \cD,$ it has to be zero on $\cD$. This is impossible because $\mu$ is a probability measure.\qed
\subsubsection*{Proof of Theorem \ref{th:sdehypo}}
We now assume that Conditions $(i),(ii), (iii)$ of Theorem \ref{th:sdehypo} hold. By Lemma~\ref{lem:hypocomp} and
Corollary~\ref{cor:main}, there exists a positive right eigenfunction for $G^\cD.$ By Conditions $(i)$ and $(iii)$ of Theorem
\ref{th:sdehypo}, there exists a point $p \in \cD$ near $\partial \cD,$ $\cD$-accessible, at which $S^j(p) \neq 0$ for some
$j \geq 1.$ Then, by Lemma~\ref{lem:hypopetiteset} and Lemma~\ref{lem:irreduc}, $(P_t^\cD)$ is irreducible. Thus, according to
Theorem~\ref{th:uniqueness}, it has a unique QSD. Such a QSD has a smooth density by Lemma \ref{lem:hyposmooth}. This proves
\bla{all the statements in} Theorem \ref{th:sdehypo}\bla{, except~\eqref{eq:cv-sde} under the additional
  assumption that there exists $x^*\in \mathcal{D}$ where the H\"ormander condition is satisfied.  This
  last result follows from Lemma \ref{lem:hyposmallset} and Theorem \ref{th:convergence}.

\section{Beyond diffusions}
\label{sec:PDMP}
This section briefly discusses a simple example of piecewise deterministic Markov process. This illustrates the applicability
of our results beyond diffusions and justifies some of the abstract conditions made in Section \ref{sec:main}. In particular, the
facts that $M$ can be chosen to be a general metric space
 rather than $\RR^n$ or sub-space of $\mathbb R^n$,
and that $D$ can be chosen to be open relative to $K = \overline{D}$ rather than open.

Let $M = \RR \times \{0,1\}$ endowed with the metric $d((x,i),(y,j))=|i-j|+|x-y|$ and  $(Z_t) = ((X_t,I_t))$  be the Markov process on $M$ defined by
$$\dot{X_t} = 2 I_t - 1,$$ and $(I_t)$ is a continuous time Markov chain on $\{0, 1\}$ having jump rates $\lambda_{0, 1} = \lambda_{1,0} = \lambda > 0.$
In words, $(X_t)$ moves at velocity $1$ either to the left or to the right and changes direction at a constant rate $\lambda.$
It is easy to check that the semi-group associated to $(Z_t)$ is Feller.

Suppose now that $(Z_t)$ is killed when $(X_t)$ exits $]0,1[.$ At the exit time, one clearly has $X_t = I_t$ almost surely, so that it is natural to set
$$D = ]0,1] \times \{0\} \cup [0,1[ \times \{1\}.$$
Defining  $K = \overline{D} = [0,1] \times \{0,1\}$, $D$ is open relative to $K$ and $$\partial_K D = \{(0,0)\} \cup \{(1,1)\}.$$
It is not hard to show that $M \setminus K$ is accessible from all $z \in K$ (in the sense of definition \ref{def:accessible}),
and that $G^D(C_b(D)) \subset C_0(D).$ The verification is left to the reader.

For any function $f : M \mapsto \RR$ and $i \in \{0,1\}$ we write $f_i(x) = f(x,i).$
\blem
\label{lem:eigenpdmp}
Let $H : [0,1] \mapsto \Rp$ and $\omega > 0$ be defined as follows:
\begin{itemize}
   \item
If $\lambda > 1, H(x) = \frac{\sin(\theta x)}{\theta}, \omega = \lambda (1 + \cos(\theta)),$
 where $0 < \theta < \pi$ is the unique solution to $\lambda \sin(\theta) = \theta;$
 \item If $\lambda = 1, H(x) = x,  \omega = 2$;
 \item If $\lambda < 1, H(x) = \frac{ \sinh(\theta x)}{\theta}, \omega = \lambda (1 + \cosh(\theta)),$ where $\theta > 0$ is the unique solution to $\lambda \sinh(\theta) = \theta.$
 \end{itemize}
 Then,
\bdes
\iti
 The map $h \in C_0(D)$ defined by
 \beq
 \label{eq:defhpdmp}
h_i(x) = (1-i) H(x) + i H(1-x)
\eeq
is a positive right eigenfunction for $G^D$ with eigenvalue $\frac{1}{\omega}.$
\itii The probability on $D$ defined by
\beq
\label{eq:defmupdmp}
\mu = \frac{1}{2 \int_0^1 H(s) ds} (h_1(x) dx \otimes \delta_0 + h_0(x) dx \otimes \delta_1)
\eeq
 is a QSD for $(Z_t)$ on $D$ with absorption rate $\omega.$
\edes
\elem
\prf
Let $L$ denote the infinitesimal generator of the process $(Z_t)$ and ${\cal D}(L)$ its domain (defined in the usual sense for Feller
processes, see e.g.~\cite{Legall2}; see also~\cite{Davis1993} more specifically on piecewise deterministic Markov processes).
 It is easy to verify that if $f :   M \mapsto \RR$ is $C^1$ in $x,$ with compact support, then $f  \in {\cal D}(L)$ and
 \beq
\label{eq:Lpdmp}
 Lf(x,i) =  \left \{ \begin{array}{c}
-f'_0(x) + \lambda (f_1(x) -f_0(x)) \mbox{ for } i = 0 \\
f'_1(x) + \lambda (f_0(x) -f_1(x)) \mbox{ for } i = 1.
 \end{array} \right.
 \eeq
Now, observe that the map $H$ satisfies the identity
\beq
\label{eq:Hidentity}
-H'(x) + \lambda (H(1-x) - H(x)) = -\omega H(x)
 \eeq  for all $0 \leq x \leq 1.$ Extending $H$ to a smooth function on $\RR$ having compact support, formula (\ref{eq:defhpdmp}) defines a map $h \in {\cal D}(L)$ such that,
 for all $z \in D$, $Lh(z) = -\omega h(z).$
Therefore, $((h(Z_{t \wedge \tau_D^{out}} ) - h(z) - \int_0^{t \wedge \tau_D^{out}} Lh(Z_s) ds))_t$ is a $\mathbb{P}_z$ martingale (because $h \in {\cal D}(L)$) bounded in $L^1,$ and by the stopping time theorem
$$- h(z) + \omega \mathbb{E}_z\left(\int_0^{\tau_D^{out}} h(Z_s) ds \right) = \mathbb{E}_z\left(h(Z_{\tau_D^{out}}) - h(z) - \int_0^{\tau_D^{out}} L h(Z_s) ds\right) = 0.$$ That is $G^D h(z) = \frac{1}{\omega} h(z).$ This prove $(i).$

We now prove $(ii).$ In order to prove that $\mu$ is a QSD with rate $\omega,$ it suffices  to show that for all $f \in C_0(D),$ $\mu G^D(f) = \frac{1}{\omega} \mu(f).$ Given such a $f,$  one can easily solve the problem $Lg = -f$ on $D$ with $g \in C_0(D).$
Namely,
$$g_0(x) =\int_0^x [f_0(u) - \lambda(x-u) (f_0(u) + f_1(u))] du + a \lambda  x$$ and
$$ g_1(x) =\int_0^x [ - f_1(u) - \lambda(x-u) (f_0(u) + f_1(u))] du + a (\lambda x + 1),$$
 where $a$ is determined by $g_1(1) = 0.$
These formulae show that $g_0$ and $g_1$ can be extended to $C^1$ maps on $\RR$ with compact support.
Reasoning exactly like in $(i)$ this shows that $G^D(f) = g$ on $D.$
The problem of showing that $\mu G^D(f) = \frac{1}{\omega} \mu(f)$ then reduces to show that $\mu (L g) = - \omega \mu(g).$ This latter identity can be easily checked,
 using (\ref{eq:defmupdmp}), (\ref{eq:Hidentity}) and integration by parts.
\qed
\brem
{\rm
 Observe  that  $(P_t^D)$ is not a semigroup on $C_b(D)$ because for any $0 < t < 1/2, \, P_t^D \Ind_D$ is discontinuous. Indeed,  $P_t^D \Ind_D(x,0) = 1$    for   $t < x < 1/2,$ while $P_t^D \Ind_D(t,0) = 1-e^{-\lambda t}.$
  It is,  however, a strongly continuous semigroup on $C_0(D).$ Let $L^D$ be its generator and ${\cal D}(L^D)$ its domain. A by product of the proof above is that ${\cal D}(L^D)$ consists of the maps $f \in C_0^1(D)$ such that
  $$f_0(0)= f_1(1) =  - f'_0(0) + \lambda f_1(0) = f'_1(1) + \lambda f_0(1) = 0,$$ and
 $L^D  f(x) = L f(x)$ for all $x\in D$ and $f \in  {\cal D}(L^D).$ Here
  $L$ is the operator defined by (\ref{eq:Lpdmp}) and  $C_0^1(D)$ stands for the set of maps $f \in C_0(D)$ such that $f_0, f_1$  can be
  extended to $C^1$ maps on $\RR.$}
\erem

\brem
{\rm
The formula defining $g = G^D(f)$ in the previous proof is also valid for $f \in C_0(D) + \RR \Ind_D.$ In particular,  $$G^D (\Ind_D)(x,0) = \mathbb{E}_{x,0}(\tau_D^{out}) = x(1 + \lambda - \lambda x).$$
Observe that the map $x \rar \mathbb{E}_{x,0}(\tau_D^{out})$ achieves its maximum at $x = 1$ for $\lambda \leq 1$ and at $x =
\frac{1+\lambda}{2\lambda}{\color{blue}<1}$ for $\lambda \geq 1$ (this last fact can be surprising at first sight).}
\erem
\bprop
\label{prop:mainpdmp}
Let $\mu$ and $h$ be as in Lemma \ref{lem:eigenpdmp}.  There exist $C, \alpha > 0$ such that for all  probability $\rho,$ on $D$
 $$\left\|\frac{\rho P_t^D}{\rho P_t^D \Ind_D} - \mu(\cdot)\right\|_{TV} \leq \frac{C}{\rho(h)} e^{-\alpha t}$$
\eprop
\prf We rely on Theorem \ref{th:convergence}. Assertions $(i)$ and $(ii)$ are satisfied. For $(iii),$ fix $\eps < 1/10$ and let $U = ]\frac{1}{2} -\eps, \frac{1}{2} + \eps [ \times \{0\}.$ Then $U$ is clearly $D$-accessible from all $z \in D.$

Set $T = \frac{1}{2} -2 \eps.$ For $z = (x,0) \in U,$ and $A$ a Borel subset of  $]\frac{1}{2} -\eps, \frac{1}{2} + \eps [,$
$$P_T^D(z, A \times \{0\}) \geq \mathbb{P}( x- \tau_1 + (\tau_2-\tau_1) - (T-\tau_2) \in A , \tau_2 < T < \tau_3),$$ where $0 <
\tau_1 < \tau_2 < \tau_3$ are the three first jump times of $I_t.$ Thus, using the fact that, given
  $\tau_2<t<\tau_3$, $(\tau_1,\tau_2)$ is distributed as the order statistics of $(U_1,U_2)$, where $U_1$ and $U_2$ are i.i.d.\
  uniform r.v.\ on $[0,t]$,
\begin{align*}
P_T^D(z, A \times \{0\}) & \geq \mathbb{P}( x- T + 2 (\tau_2-\tau_1)  \in A | \tau_2 < T < \tau_3) \frac{(\lambda T)^2}{2} e^{-\lambda T} \\ &
= \frac{(\lambda T)^2}{2} e^{-\lambda T}  \int_A \frac{(x + T -u)}{2T^2} du  \geq c Leb(A)
\end{align*}
for some $c > 0.$ This shows condition $(iii)$ of Theorem \ref{th:convergence} with $d \xi = c \Ind_U du \otimes \delta_0.$
\qed
\brem {\em The minorization of $P_t^D(z,\cdot)$ proved in Proposition \ref{prop:mainpdmp} could also be obtained for  general PDMPs
  under appropriate bracket conditions  on the vector fields defining the PDMP, along the lines of the general results proved in
  \cite{BakhtinHurth} or \cite{BLMZ15}. Hence Theorem~\ref{th:convergence} applies to such processes provided one can
    prove that $G^D(C_b(D))\subset C_0(D)$ and $G^D$ is a compact operator on $C_0(D)$, by Corollary~\ref{cor:main}.}
\erem

\subsection*{Acknowledgment}
We thank Emmanuel Trelat for fruitful discussions and his suggestion to use reference \cite{rothschild1976}, in Lemma \ref{lem:hypocomp}.
We thank Josef Hofbauer for his valuable suggestions on Lemma \ref{lem:attractor}. We also thank 4 anonymous referees for their careful reading and valuable suggestions to improve  this paper.
This research is supported by  the Swiss National Foundation  grants 200020 196999 and 200020 219913.
\bibliography{Killedprocesses}
\bibliographystyle{amsplain}

\end{document}